# Accelerated Stochastic Greedy Coordinate Descent by Soft Thresholding Projection onto Simplex

Chaobing Song, Shaobo Cui, Shu-Tao Xia, Yong Jiang

**Abstract**

In this paper we study the well-known greedy coordinate descent (GCD) algorithm to solve $\ell_1$-regularized problems and improve GCD by the two popular strategies: Nesterov's acceleration and stochastic optimization. Firstly, we propose a new rule for greedy selection based on an $\ell_1$-norm square approximation which is nontrivial to solve but convex; then an efficient algorithm called "SOft ThreshOlding PrOjection (SOTOPO)" is proposed to exactly solve the $\ell_1$-regularized $\ell_1$-norm square approximation problem, which is induced by the new rule. Based on the new rule and the SOTOPO algorithm, the Nesterov's acceleration and stochastic optimization strategies are then successfully applied to the GCD algorithm. The resulted algorithm called accelerated stochastic greedy coordinate descent (ASGCD) has the optimal convergence rate $O(\sqrt{1/\epsilon})$; meanwhile, it reduces the iteration complexity of greedy selection up to a factor of sample size. Both theoretically and empirically, we show that ASGCD has better performance for high-dimensional and dense problems with sparse solution.

## I. INTRODUCTION

In large-scale convex optimization, first-order methods are widely used due to their cheap iteration cost. In order to improve the convergence rate and reduce the iteration cost further, two important strategies are used in first-order methods: Nesterov's acceleration and stochastic optimization. Nesterov's acceleration schemes are referred to techniques that uses some algebra trick to accelerate first-order algorithms; while stochastic optimization is referred to optimization methods that samples one training example or one dual coordinate at random from the training data in each iteration. Assume the objective function $F(x)$ is convex and smooth. Let $F^* = \min_{x \in R^d} F(x)$ be the optimal value. In order to find an approximate solution $x$ that satisfies $F(x) - F^* \leq \epsilon$, the vanilla gradient descent method needs $O(1/\epsilon)$ iterations. While after applying the Nesterov's acceleration scheme [15], the resulted accelerated full gradient method (AFG) [15] only needs $O(\sqrt{1/\epsilon})$ iterations, which is optimal for first-order algorithms [15]. Meanwhile, assume $F(x)$ is also a finite sum of $n$ sample convex functions. By sampling one training example, the resulted stochastic gradient descent (SGD) and its variants [12, 21, 1] can reduce the iteration complexity by a factor of the sample size. As an alternative of SGD, randomized coordinate descent (RCD) can also reduces the iteration complexity by a factor of the sample size [14] and obtain the optimal convergence rate $O(\sqrt{1/\epsilon})$ by Nesterov's acceleration [13, 22]. The development of gradient descent and RCD raises an interesting problem: can the Nesterov's acceleration and stochastic optimization strategies be used to improve other existing first-order algorithms?

In this paper, we answer this question partly by studying coordinate descent with Gauss-Southwell selection, i.e., greedy coordinate descent (GCD). GCD is widely used for solving sparse optimization problem in machine learning [20, 9, 16]. If the optimization problem has a sparse solution, it is more suitable than its counterpart RCD. However, the theoretical convergence rate is still $O(1/\epsilon)$. Meanwhile if the iteration complexity is comparable, GCD will be preferable than RCD [16]. However in the general case, in order to do exact Gauss-Southwell selection, computing the full gradient beforehand is necessary, which causes GCD has much higher iteration complexity than RCD. To be concrete, in this paper we consider the well-known nonsmooth $\ell_1$-regularized problem:

$$\min_{x \in \mathbb{R}^d} \left\{ F(x) \stackrel{\text{def}}{=} f(x) + \lambda \|x\|_1 \stackrel{\text{def}}{=} \frac{1}{n} \sum_{j=1}^n f_j(x) + \lambda \|x\|_1 \right\}, \qquad (1)$$

where $\lambda \geq 0$ is a regularization parameter, $f(x) = \frac{1}{n} \sum_{j=1}^n f_j(x)$ is a smooth convex function that is a finite average of $n$ smooth convex function $f_j(x)$. Given samples $\{(a_1, b_1), (a_2, b_2), \ldots, (a_n, b_n)\}$ with $a_j \in \mathbb{R}^d, b_j \in \mathbb{R}$ ($j \in [n] \stackrel{\text{def}}{=} \{1, 2, \ldots, n\}$), if each $f_j(x) = f_j(a_j^T x, b_j)$, then (1) is an $\ell_1$-regularized empirical risk minimization ($\ell_1$-ERM) problem. If $f_j(x) = \frac{1}{2}(b_j - a_j^T x)^2$, then (1) is Lasso; if $f_j(x) = \log(1 + \exp(-b_j a_j^T x))$, then $\ell_1$-regularized logistic regression is obtained.

In the above nonsmooth case, the Gauss-Southwell rule has 3 different variants [16, 20]: GS-$s$, GS-$r$ and GS-$q$. The GCD algorithm with all the 3 rules can be viewed as the following procedure: in each iteration based on a quadratic approximation of $f(x)$ in (1), minimizing a surrogate objective function under the constraint that the direction vector used for update has at most 1 nonzero entries. The resulted problems under the 3 rules are easy to solve but are *nonconvex* due to the cardinality constraint of direction vector. While when using Nesterov's acceleration scheme, convexity is needed for the derivation of the optimal convergence rate $O(\sqrt{1/\epsilon})$ [15]. Therefore, it is impossible to accelerate GCD by the Nesterov's acceleration scheme under the existing 3 rules.

In this paper, we propose a novel variant of Gauss-Southwell rule by using an $\ell_1$-norm sqoximation of $f(x)$ rather than quadratic approximation. The resulted $\ell_1$-*regularized* $\ell_1$-*norm square approximation* problem is nontrivial to solve but is *convex*. The main challenge in this paper is to solve the nontrivial $\ell_1$-regularized $\ell_1$-norm square approximation problem. In this paper, we propose an efficient SOft PrOjection PrOjection (SOTOPO) algorithm to exactly solve this problem. The SOTOPO algorithm has $O(d + |Q| \log |Q|)$, where it is often the case $|Q| \ll d$. The complexity result $O(d + |Q| \log |Q|)$ is better than $O(d \log d)$ of its counterpart SOPOPO [17], which is an Euclidean projection method on polyhedra.

Then based on this new rule and the SOTOPO algorithm, we accelerate GCD to attain the optimal convergence rate $O(\sqrt{1/\epsilon})$ by combing a delicately selected mirror descent step. Meanwhile, we show that it is not necessary to compute full gradient beforehand: sampling one training example and computing a noisy gradient rather than full gradient is enough to perform greedy selection. This stochastic optimization technique reduces the iteration complexity of greedy selection by a factor of the sample size. The final result is an accelerated stochastic greedy coordinate descent (ASGCD) algorithm.

Assume $x^*$ is an optimal solution of (1). Assume that $f_j(x)$(for all $j \in [n]$) are $L_p$-smooth w.r.t. $\|\cdot\|_p$ $(p = 1, 2)$, i.e., for all $x, y \in \mathbb{R}^d$,
$$\|\nabla f_j(x) - \nabla f_j(y)\|_q \leq L_p \|x - y\|_p, \tag{2}$$
where if $p = 1$, then $q = \infty$; if $p = 2$, then $q = 2$.

In order to find an $x$ that satisfies $F(x) - F(x^*) \leq \epsilon$, ASGCD needs $O\left(\frac{\sqrt{CL_1}\|x^*\|_1}{\sqrt{\epsilon}}\right)$ iterations, where $C$ is a function of $d$ that varies slowly over $d$ and is upper bounded by $\log^2(d)$. For high-dimensional and dense problems with sparse solution, ASGCD has better performance than the state of the art. Experiments demonstrates the theoretical result.

Notations: Let $[d]$ denote the set $\{1, 2, \ldots, d\}$. Let $\mathbb{R}_+$ denote the set of nonnegative real number. For $x \in \mathbb{R}^d$, let $\|x\|_p = (\sum_{i=1}^d |x_i|^p)^{\frac{1}{p}}$ $(1 \leq p < \infty)$ denote the $\ell_p$-norm and $\|x\|_\infty = \max_{i \in [d]} |x_i|$ denote the $\ell_\infty$-norm of $x$. For a vector $x$, let $\dim(x)$ denote the dimension of $x$; let $x_i$ denote the $i$-th element of $x$. For a gradient vector $\nabla f(x)$, let $\nabla_i f(x)$ denote the $i$-th element of $\nabla f(x)$. For a set $S$, let $|S|$ denote the cardinality of $S$. Denote the simplex $\triangle_d = \{\theta \in \mathbb{R}_+^d : \sum_{i=1}^d \theta_i = 1\}$.

## II. THE SOTOPO ALGORITHM

The proposed SOTOPO algorithm aims to solve the proposed new rule, *i.e.*, minimizing an $\ell_1$-regularized $\ell_1$-norm square approximation problem, which is described as the following iteration,
$$\tilde{h} \stackrel{\text{def}}{=} \arg\min_{g \in \mathbb{R}^d} \left\{ \langle \nabla f(x), g \rangle + \frac{1}{2\eta} \|g\|_1^2 + \lambda \|x + g\|_1 \right\}, \tag{3}$$
$$\tilde{x} \stackrel{\text{def}}{=} x + \tilde{h}, \tag{4}$$
where $x$ denotes the current iteration, $\eta$ denotes a step size, $\tilde{h}$ denotes the director vector for update and $\tilde{x}$ denotes the next iteration. The number of nonzero entries of $\tilde{h}$ denotes how many coordinates will be updated in this iteration. Unlike the quadratic approximation used in GS-$s$, GS-$r$ and GS-$q$ rules, in the new rule the coordinate(s) to update is implicitly selected by the sparsity-inducing property of the $\ell_1$-norm square $\|g\|_1^2$ rather than using the cardinality constraint $\|g\|_0 \leq 1$ [16, 20]. By [7, §9.4.2], when the nonsmooth term $\lambda\|x + g\|_1$ in (1) does not exist, the minimizer of the $\ell_1$-norm square approximation (i.e., $\ell_1$-norm steepest descent) is equivalent to GCD. When $\lambda\|x + g\|_1$ exists, generally, there may be one or more coordinates to update in this new rule. In addition, (3) is an unconstrained problem and thus is feasible.

**Remark 1.** *When $\tilde{h}$ is not unique, we can choose one solution arbitrarily without influencing the theoretical analysis; as an alternative, we can also choose the one with the minimal nonzero entries so as to update the coordinates of $x$ as small as possible in each iteration. For simplicity, in this paper, we aim to find an arbitrary minimizer of* (3).

### A. A variational reformulation and its properties

(3) involves the nonseparable and nonsmooth term $\|g\|_1^2$ and the nonsmooth term $\|x + g\|_1$. While by the variational identity $\|g\|_1^2 = \inf_{\theta \in \triangle_d} \sum_{i=1}^d \frac{g_i^2}{\theta_i}$ in [4] [1], we give Lemma 1.

**Lemma 1.** *By defining*
$$J(g, \theta) \stackrel{\text{def}}{=} \langle \nabla f(x), g \rangle + \frac{1}{2\eta} \sum_{i=1}^d \frac{g_i^2}{\theta_i} + \lambda\|x + g\|_1, \tag{5}$$
$$\tilde{g}(\theta) \stackrel{\text{def}}{=} \arg\min_{g \in \mathbb{R}^d} J(g, \theta), \quad J(\theta) \stackrel{\text{def}}{=} J(\tilde{g}(\theta), \theta), \tag{6}$$
$$\tilde{\theta} \stackrel{\text{def}}{=} \arg\inf_{\theta \in \triangle_d} J(\theta), \tag{7}$$

---
[1]The infima can be replaced by minimization if the convention "0/0 = 0" is used.

then the minimization problem of $\tilde{h}$ in (3) is equivalent to the problem (7) with $\tilde{h} = \tilde{g}(\tilde{\theta})$. Meanwhile, the $\tilde{g}(\theta)$ and $J(\theta)$ in (6) are both coordinate separable with the expressions

$$\forall i \in [d], \quad \tilde{g}_i(\theta) = \tilde{g}_i(\theta_i) \stackrel{\text{def}}{=} \text{sign}(x_i - \theta_i \eta \nabla_i f(x)) \cdot \max\{0, |x_i - \theta_i \eta \nabla_i f(x)| - \theta_i \eta \lambda\} - x_i, \tag{8}$$

$$J(\theta) = \sum_{i=1}^{d} J_i(\theta_i), \quad \text{where} \quad J_i(\theta_i) \stackrel{\text{def}}{=} J(\tilde{g}_i(\theta_i), \theta_i). \tag{9}$$

In Lemma 1, (8) is obtained by the iterative soft thresholding operator [6]. By Lemma 1, we can reformulate (3) into the problem (5) about two parameter $g$ and $\theta$. Then by the joint convexity, we swap the optimization order of $g$ and $\theta$. Fixing $\theta$ and optimizing with respect to $g$, we can get a closed form of $\tilde{g}(\theta)$, which is a vector function about $\theta$. Substituting $\tilde{g}(\theta)$ into $J(g, \theta)$, we get the problem (7) about $\theta$. Finally, the optimal solution of $\tilde{h}$ in (3) can be obtained by $\tilde{h} = \tilde{g}(\tilde{\theta})$.

The explicit expression of each $J_i(\theta_i)$ can be given by substituting (8) into (9). Because $\theta \in \triangle_d$, we have for all $i \in [d]$, $0 \leq \theta_i \leq 1$. To deduce the SOTOPO algorithm, it is observed that the derivate $J'_i(\theta_i)$ has the following properties.

**Lemma 2.** *Assume that for all $i \in [d]$, $J'_i(0)$ and $J'_i(1)$ have been computed. Denote $r_{i1} \stackrel{\text{def}}{=} \frac{|x_i|}{\sqrt{-2\eta J'_i(0)}}$ and $r_{i2} \stackrel{\text{def}}{=} \frac{|x_i|}{\sqrt{-2\eta J'_i(1)}}$, then $J'_i(\theta_i)$ belongs to one of the 4 cases,*

$$(\text{case } a): J'_i(\theta_i) = 0, \quad 0 \leq \theta_i \leq 1, \qquad (\text{case } b): J'_i(\theta_i) = J'_i(0) < 0, \quad 0 \leq \theta_i \leq 1,$$

$$(\text{case } c): J'_i(\theta_i) = \begin{cases} J'_i(0), & 0 \leq \theta_i \leq r_{i1} \\ -\frac{x_i^2}{2\eta\theta_i^2}, & r_{i1} < \theta_i \leq 1 \end{cases}, \quad (\text{case } d): J'_i(\theta_i) = \begin{cases} J'_i(0), & 0 \leq \theta_i \leq r_{i1} \\ -\frac{x_i^2}{2\eta\theta_i^2}, & r_{i1} < \theta_i < r_2 \\ J'_i(1), & r_{i2} \leq \theta_i \leq 1 \end{cases}.$$

Lemma 2 shows that $J'_i(\theta_i)$ can be a constant function or a piecewise function. Although the formulation of $J'_i(\theta_i)$ is complicated, by summarizing the property of the 4 cases in Lemma 2, we have Corollary 1.

**Corollary 1.** *For all $i \in [d]$ and $0 \leq \theta_i \leq 1$, if the derivate $J'_i(\theta_i)$ is not always 0, then $J'_i(\theta_i)$ is a non-decreasing, continuous function with value always less than 0.*

### B. The property of the optimal solution

The Lagrangian of the problem (7) is

$$\mathcal{L}(\theta, \gamma, \zeta) \stackrel{\text{def}}{=} J(\theta) + \gamma \left( \sum_{i=1}^{d} \theta_i - 1 \right) - \langle \zeta, \theta \rangle, \tag{10}$$

where $\gamma \in \mathbb{R}$ is a Lagrange multiplier and $\zeta \in \mathbb{R}_+^d$ is a vector of non-negative Lagrange multipliers. Due to coordinate separable property of $J(\theta)$ in (9), it follows that $\frac{\partial J(\theta)}{\partial \theta_i} = J'_i(\theta_i)$. Then the KKT condition of (10) can be written as

$$\forall i \in [d], \quad J'_i(\theta_i) + \gamma - \zeta_i = 0, \quad \zeta_i \theta_i = 0, \quad \text{and} \quad \sum_{i=1}^{d} \theta_i = 1. \tag{11}$$

Reformulate the KKT condition (11), we have Lemma 3.

**Lemma 3.** *If $(\tilde{\gamma}, \tilde{\theta}, \tilde{\zeta})$ is a stationary point of the problem (10), then $\tilde{\theta}$ is an optimal solution of (7). Meanwhile, denote $S \stackrel{\text{def}}{=} \{i : \tilde{\theta}_i > 0\}$ and $T \stackrel{\text{def}}{=} \{j : \tilde{\theta}_j = 0\}$, then the KKT condition can be formulated as*

$$\begin{cases} \sum_{i \in S} \tilde{\theta}_i = 1; \\ \text{for all } j \in T, \quad \tilde{\theta}_j = 0; \\ \text{for all } i \in S, \quad \tilde{\gamma} = -J'_i(\tilde{\theta}_i) \geq \max_{j \in T} -J'_j(0). \end{cases} \tag{12}$$

### C. The soft thresholding projection algorithm

In Lemma 3, it is shown that the negative derivates $-J'_i(\tilde{\theta}_i)$ with $\tilde{\theta}_i > 0$ are equal to a single variable $\tilde{\gamma}$. Therefore, a much simpler problem can be obtained if we know the coordinates of these positive elements. At first glance, it seems difficult to identify these coordinates, because the number of potential subsets of coordinates is clearly exponential on the dimension $d$. However, the property clarified by Lemma 2 enables an efficient procedure for identifying the non-zero elements of $\tilde{\theta}$. Lemma 4 is a key tool in deriving the procedure for identifying the non-zero elements of $\tilde{\theta}$.

**Lemma 4** (Non-zero element identification). *Let $\tilde{\theta}$ be an optimal solution of (7). Let $s$ and $t$ be two coordinates such that $J'_s(0) < J'_t(0)$. If $\tilde{\theta}_s = 0$, then $\tilde{\theta}_t$ must be 0 as well; equivalently, if $\tilde{\theta}_t > 0$, then $\tilde{\theta}_s$ must be greater than 0 as well.*

Lemma 4 shows that if we sort $u \stackrel{\text{def}}{=} -\nabla J(0)$ such that $u_{i_1} \geq u_{i_2} \geq \cdots \geq u_{i_d}$, where $\{i_1, i_2, \ldots, i_d\}$ is a permutation of $[d]$, then the set $S$ in Lemma 3 is of the form $\{i_1, i_2, \ldots, i_\varrho\}$, where $1 \leq \varrho \leq d$. If $\varrho$ is obtained, then we can use the fact that for all $j \in [\varrho]$,

$$-J'_{i_j}(\tilde{\theta}_{i_j}) = \tilde{\gamma} \quad \text{and} \quad \sum_{j=1}^{\varrho} \tilde{\theta}_{i_j} = 1 \tag{13}$$

to compute $\tilde{\gamma}$. Therefore, by Lemma 4, we can efficiently identify the nonzero elements of the optimal solution $\tilde{\theta}$ after a sort operation, which costs $O(d \log d)$. While based on the property of $J'_i(\theta_i)$ in Lemma 2 and the KKT condition in Lemma 2, the sort cost $O(d \log d)$ can be further reduced by the following Lemma 5.

**Lemma 5** (Efficient identification). *Assume $\tilde{\theta}$ and $S$ are given in Lemma 3. Then for all $i \in S$,*

$$-J'_i(0) \geq \max_{j \in [d]} \{-J'_j(1)\}. \tag{14}$$

By Lemma 5, before ordering $u$, we can filter out all the coordinates $i$'s that satisfy $-J'_i(0) < \max_{j \in [d]} -J'_j(1)$. Based on Lemmas 4 and 5, we propose the SOft ThreshOlding PrOjection (SOTOPO) algorithm in Alg. 1 to efficiently obtain an optimal solution $\tilde{\theta}$. In the step 1, by Lemma 5, we find the quantity $v_m, i_m$ and $Q$. In the step 2, by Lemma 4, we sort the elements $\{-J'_i(0)| i \in Q\}$. In the step 3, because $S$ in Lemma 3 is of the form $\{i_1, i_2, \ldots, i_\varrho\}$, we search the quantity $\rho$ from 1 to $|Q| + 1$ until a stopping criteria is met. In Alg. 1, $\rho$ or $\rho - 1$ may be the number of non-zero elements of $\tilde{\theta}$. In the step 4, we compute the $\tilde{\gamma}$ in Lemma 3 according to the conditions. In the step 5, the optimal $\tilde{\theta}$ and the corresponding $\tilde{h}, \tilde{x}$ are given.

---

**Algorithm 1** $\tilde{x} = \text{SOTOPO}(\nabla f(x), x, \lambda, \eta)$

1) Find

$$(v_m, i_m) \stackrel{\text{def}}{=} (\max_{i \in [d]} \{-J'_i(1)\}, \arg\max_{i \in [d]} \{-J'_i(1)\}), \quad Q \stackrel{\text{def}}{=} \{i \in [d] | -J'_i(0) > v_m\}.$$

2) Sort $\{-J'_i(0)| i \in Q\}$ such that $-J'_{i_1}(0) \geq -J'_{i_2}(0) \geq \cdots \geq -J'_{i_{|Q|}}(0)$, where $\{i_1, i_2, \ldots, i_{|Q|}\}$ is a permutation of the elements in $Q$. Denote

$$v \stackrel{\text{def}}{=} (-J'_{i_1}(0), -J'_{i_2}(0), \ldots, -J'_{i_{|Q|}}(0), v_m), \quad \text{and set } i_{|Q|+1} \stackrel{\text{def}}{=} i_m, v_{|Q|+1} = v_m.$$

3) For $j \in [|Q| + 1]$, denote $R_j = \{i_k | k \in [j]\}$. Search from 1 to $|Q| + 1$ to find the quantity

$$\rho \stackrel{\text{def}}{=} \min\left\{j \in [|Q| + 1]| J'_{i_j}(0) = J'_{i_j}(1) \text{ or } \sum\nolimits_{l \in R_j} |x_l| \geq \sqrt{2\eta v_j} \text{ or } j = |Q| + 1\right\}.$$

4) The $\tilde{\gamma}$ in Lemma 3 are given by

$$\tilde{\gamma} = \begin{cases} \left(\sum_{l \in R_{\rho-1}} |x_l|\right)^2/(2\eta), & \text{if } \sum_{l \in R_{\rho-1}} |x_l| \geq \sqrt{2\eta v_\rho}; \\ v_\rho, & \text{otherwise.} \end{cases}$$

5) Then the $\tilde{\theta}$ in Lemma 3 and its corresponding $\tilde{h}, \tilde{x}$ in (3) and (4) are obtained by

$$(\tilde{\theta}_l, \tilde{h}_l, \tilde{x}_l) = \begin{cases} \left(\frac{|x_l|}{\sqrt{2\eta\tilde{\gamma}}}, x_l, 0\right), & \text{if } l \in R_\rho \setminus \{i_\rho\}; \\ \left(1 - \sum_{k \in R_\rho \setminus \{i_\rho\}} \tilde{\theta}_k, \tilde{g}_l(\tilde{\theta}_l), x_l + \tilde{g}_l(\tilde{\theta}_l)\right), & \text{if } l = i_\rho; \\ (0, 0, x_l), & \text{if } l \in [d] \setminus R_\rho. \end{cases}$$

---

In Theorem 1, we give the main result about the SOTOPO algorithm.

**Theorem 1.** *The SOTOPO algorihtm in Alg. 1 can get the exact minimizer $\tilde{h}, \tilde{x}$ of the $\ell_1$-regularized $\ell_1$-norm square approximation problem in (3) and (4).*

The SOTOPO algorithm seems complicated but is indeed efficient. The dominant operations in Alg. 1 are steps 1 and 2 with the total cost $O(d + |Q| \log |Q|)$. To show the effect of the complexity reduction by Lemma 5, we give the following fact.

**Proposition 1.** *For the optimization problem defined in (5)-(7), where $\lambda$ is the regularization parameter of the original problem (1), we have that*

$$0 \leq \max_{i \in [d]} \{-J'_i(0)\} - \max_{j \in [d]} \{-J'_j(1)\} \leq 2\lambda. \tag{15}$$

By Lemma 5 and Proposition 1, we have that for all $i \in Q$, $\max_{k \in [d]} \{-J'_k(0)\} \leq -J'_i(0) + 2\lambda$. Therefore after applying Lemma 5, in the step 1 of Alg. 1, the maximal difference between two elements in $\{-J'_i(0)| i \in Q\}$ in the step 2 of Alg. 1

---

**Algorithm 2** $(\tilde{x}, \tilde{\vartheta}) = \text{pCOMID}(g, \vartheta, q, \lambda, \alpha)$

1) $\forall i \in [d], \tilde{\vartheta}_i = \text{sign}(\vartheta_i - \alpha g_i) \cdot \max\{0, |\vartheta_i - \alpha g_i| - \alpha\lambda\}$;
2) $\forall i \in [d], \tilde{x}_i = \frac{\text{sign}(\tilde{\vartheta}_i)|\tilde{\vartheta}_i|^{q-1}}{\|\tilde{\vartheta}\|_q^{q-2}}$;
3) **Output:** $\tilde{x}, \tilde{\vartheta}$.

---

will be less than $2\lambda$. Therefore at least the coordinates $j$'s that satisfy $\max_{k \in [d]}\{-J'_k(0)\} > -J'_j(0) + 2\lambda$ will be not contained in $Q$. In practice, it can considerably reduce the sort complexity.

**Remark 2.** *SOTOPO can be viewed as an extension of the SOPOPO algorithm [17] by changing the objective function from Euclidean distance to a more general function $J(\theta)$ in (9). It should be noted that Lemma 5 does not have counterpart in the case that the objective function is Euclidean distance [17]. In addition, some extension of randomized median finding algorithm [10] with linear time in our setting is also deserved to research, due to the limited space, it is left for further discussion.*

## III. THE ASGCD ALGORITHM

---

**Algorithm 3** ASGCD

$\delta = \log(d) - 1 - \sqrt{(\log(d) - 1)^2 - 1}$;
$p = 1 + \delta, q = \frac{p}{p-1}, C = \frac{d^{\frac{2\delta}{1+\delta}}}{\delta}$;
$z_0 = y_0 = \tilde{x}_0 = \vartheta_0 = 0$;
$\tau_2 = \frac{1}{2}, m = \lceil \frac{n}{b} \rceil, \eta = \frac{1}{(1 + 2\frac{n-b}{b(n-1)})L_1}$;
**for** $s = 0, 1, 2, \ldots, S-1$, do

1) $\tau_{1,s} = \frac{2}{s+4}, \alpha_s = \frac{\eta}{\tau_{1,s}C}$;
2) $\mu_s = \nabla f(\tilde{x}_s)$;
3) **for** $l = 0, 1, \ldots, m-1$, do
   a) $k = (sm) + l$;
   b) randomly sample a mini batch $\mathcal{B}$ of size b from $\{1, 2, \ldots, n\}$ with equal probability;
   c) $x_{k+1} = \tau_{1,s}z_k + \tau_2\tilde{x}_s + (1 - \tau_{1,s} - \tau_2)y_k$;
   d) $\tilde{\nabla}_{k+1} = \mu_s + \frac{1}{b}\sum_{j \in \mathcal{B}}(\nabla f_j(x_{k+1}) - \nabla f_j(\tilde{x}_s))$;
   e) $y_{k+1} = \text{SOTOPO}(\tilde{\nabla}_{k+1}, x_{k+1}, \lambda, \eta)$;
   f) $(z_{k+1}, \vartheta_{k+1}) = \text{pCOMID}(\tilde{\nabla}_{k+1}, \vartheta_k, q, \lambda, \alpha_s)$;
   **end for**
4) $\tilde{x}_{s+1} = \frac{1}{m}\sum_{l=1}^{m} y_{sm+l}$;

**end for**
**Output:** $\tilde{x}_S$

---

Now we can come back to our motivation, *i.e.*, accelerate GCD to obtain the optimal convergence rate $O(1/\sqrt{\epsilon})$ by Nesterov's acceleration and reduce the complexity of greedy selection by stochastic optimization. The main idea is that although like any (block) coordinate descent algorithm, the proposed new rule, *i.e.*, minimizing the problem in (3), performs update on one or several coordinates, it is a generalized proximal gradient descent problem based on $\ell_1$-norm. Therefore this rule can be applied into the existing Nesterov's acceleration and stochastic optimization framework "Katyusha" [1] if it can be solved efficiently. The final result is the accelerated stochastic greedy coordinate descent (ASGCD) algorithm, which is described in Alg. III.

In Alg. III, the gradient descent step $3(e)$ is solved by the proposed SOTOPO algorithm, while the mirror descent step $3(f)$ is solved by the COMID algorithm with $p$-norm divergence [11, Sec. 7.2]. We denote the mirror descent step as pCOMID in Alg. 2. All other parts are standard steps in the Katyusha framework except some parameter settings. For example, instead of the custom setting $p = 1 + 1/\log(d)$ [18, 11], a particular choice $p = 1 + \delta$ ($\delta$ is defined in Alg. III) is used to minimize $C = \frac{d^{\frac{2\delta}{1+\delta}}}{\delta}$. $C$ varies slowly over $d$ and is upper bounded by $\log^2(d)$. Meanwhile, $\alpha_{k+1}$ depends on the extra constant $C$. Furthermore, the step size $\eta = \frac{1}{(1 + 2\frac{n-b}{b(n-1)})L_1}$ is used, where $L_1$ is defined in (2). Finally, unlike [1, Alg. 2], we let the batch size $b$ as an algorithm parameter to cover both the stochastic case $b < n$ and the deterministic case $b = n$. To the best of our knowledge, the existing GCD algorithms are deterministic, therefore by setting $b = n$, we can compare with the existing GCD algorithms better.

Based on the efficient SOTOPO algorithm, ASGCD has nearly the same iteration complexity with the normal form [1, Alg. 2] of Katyusha. Meanwhile we have the following convergence rate.

**Theorem 2.** *If each $f_j(x)(j \in [n])$ is convex, $L_1$-smooth in (2) and $x^*$ is an optimum of the $\ell_1$-regularized problem (1), then ASGCD satisfies*

$$\mathbb{E}[F(\tilde{x}^S)] - F(x^*) \leq \frac{4}{(S+3)^2}\left(1 + \frac{1+2\beta(b)}{2m}C\right)L_1\|x^*\|_1^2 = O\left(\frac{CL_1\|x^*\|_1^2}{S^2}\right), \quad (16)$$

*where $\beta(b) = \frac{n-b}{b(n-1)}$, $S$, $b$, $m$ and $C$ are given in Alg. III. In other words, ASGCD achieves an $\epsilon$-additive error (i.e., $\mathbb{E}[F(\tilde{x}^S)] - F(x^*) \leq \epsilon$) using at most $O\left(\frac{\sqrt{CL_1}\|x^*\|_1}{\sqrt{\epsilon}}\right)$ iterations.*

The bound depends on $\|x^*\|_1$ rather than $\|x^*\|_2$ and on $L_1$ rather than $L_2$ ($L_1$ and $L_2$ are defined in (2)). To the best of our knowledge, for the $\ell_1$-regularized problem in (1), the state of the art convergence rate of first-order algorithms is $O\left(\frac{\sqrt{L_2}\|x^*\|_2}{\sqrt{\epsilon}}\right)$, obtained by [1, Alg. 2]. For the $\ell_1$-ERM problem, if the samples are high-dimensional and dense and the regularization parameter $\lambda$ is relatively large, then it is possible that $L_1 \ll L_2$ (in the extreme case, $L_2 = dL_1$ [9]) and $\|x^*\|_1 \approx \|x^*\|_2$. It means that ASGCD will have better theoretical guarantee in this case.

**Remark 3.** *When the batch size $b = n$, ASGCD is a deterministic algorithm. In this case, we can use a better smooth constant $T_1$ that satisfies $\|\nabla f(x) - \nabla f(y)\|_\infty \leq T_1\|x-y\|_1$ rather than $L_1$.*

**Remark 4.** *The necessity of computing the full gradient beforehand is the main bottleneck of GCD in application [16]. There exists some work [9] to avoid the computation of full gradient by performing some approximate greedy selection. While the method in [9] needs preprocessing, incoherence condition for dataset and somewhat complicated. Contrary to [9], the proposed ASGCD algorithm reduce the complexity of greedy selection by a factor up to $n$ in terms of the amortized cost by simply applying the existing stochastic variance reduction framework.*

## IV. EXPERIMENTS

In this section, we use numerical experiments to demonstrate the theoretical results in Section III and show the empirical performance of ASGCD with batch size $b=1$ (In Fig. 1 it is denoted as ASGCD ($b=1$)) and its deterministic version with $b=n$ (In Fig. 1 it is denoted as ASGCD ($b=n$)) respectively. In addition, following the claims to using data access rather than CPU time [18] and the recent SGD and RCD literature [22, 12, 13, 1], we use the data access, i.e., the number of times the algorithm accesses the data matrix, to measure the algorithm performance. To show the effect of Nesterov's acceleration, we compare ASGCD ($b=n$) with the non-accelerated greedy coordinate descent with GS-$q$ rule, i.e., coordinate gradient descent (CGD) [20]. To show the effect of both Nesterov's acceleration and stochastic strategy, we compare ASGCD($b=1$) with Katyusha [1, Alg. 2]. To show the effect of the proposed new rule in Section II that based on $\ell_1$-norm square approximation, we compare ASGCD ($b=n$) with the $\ell_2$-norm based proximal accelerated full gradient (AFG) implemented by the linear coupling framework [3]. Meanwhile, as a benchmark of stochastic optimization for the problem with finite-sum structure, we also show the performance of proximal stochastic variance reduced gradient (SVRG) [21]. In addition, based on the [1] and our experiments, we find that "Katyusha" [1, Alg. 2] has the best empirical performance in general for the $\ell_1$-regularized problem (1). Therefore other well-known state-of-art algorithms, such as APCG [13] and accelerated SDCA [19], are not included in the experiments.

The datasets are obtained from LIBSVM data [8] and summarized in Table I. All the algorithms are used to solve the following lasso problem

$$\min_{x \in \mathbb{R}^d}\{f(x) + \lambda\|x\|_1 = \frac{1}{2n}\|b - Ax\|_2^2 + \lambda\|x\|_1\} \quad (17)$$

on the 3 datasets, where $A = (a_1, a_2, \ldots, a_n)^T = (h_1, h_2, \ldots, h_d) \in \mathbb{R}^{n \times d}$ with each $a_j \in \mathbb{R}^d$ represents a sample vector and $h_i \in \mathbb{R}^n$ represents a feature vector, $b \in \mathbb{R}^n$ is the prediction vector.

For ASGCD ($b=1$) and Katyusha [1, Alg. 2], we can use the tight smooth constant $L_1 = \max_{j \in [n], i \in [d]} |a_{j,i}^2|$, $L_2 = \max_{j \in [n]} \|a_j\|_2^2$ respectively in their implementation. While for ASGCD ($b=n$) and AFG, the better smooth constant $T_1 = \frac{\max_{i \in [d]} \|h_i\|_2^2}{n}$, $T_2 = \frac{\|A\|^2}{n}$, are used respectively. The learning rate of CGD and SVRG are tuned in $\{10^{-6}, 10^{-5}, 10^{-4}, 10^{-3}, 10^{-2}, 10^{-1}\}$.

We use $\lambda = 10^{-6}$ and $\lambda = 10^{-2}$ in the experiments. In addition, for each case (Dataset, $\lambda$), AFG is used to find an optimum $x^*$ with enough accuracy.

TABLE I
CHARACTERISTICS OF THREE REAL DATASETS.

| Dataset Name | # Samples $n$ | # Features $d$ |
|---|---|---|
| Leukemia | 38 | 7129 |
| Gisette | 6000 | 5000 |
| Mnist | 60000 | 780 |

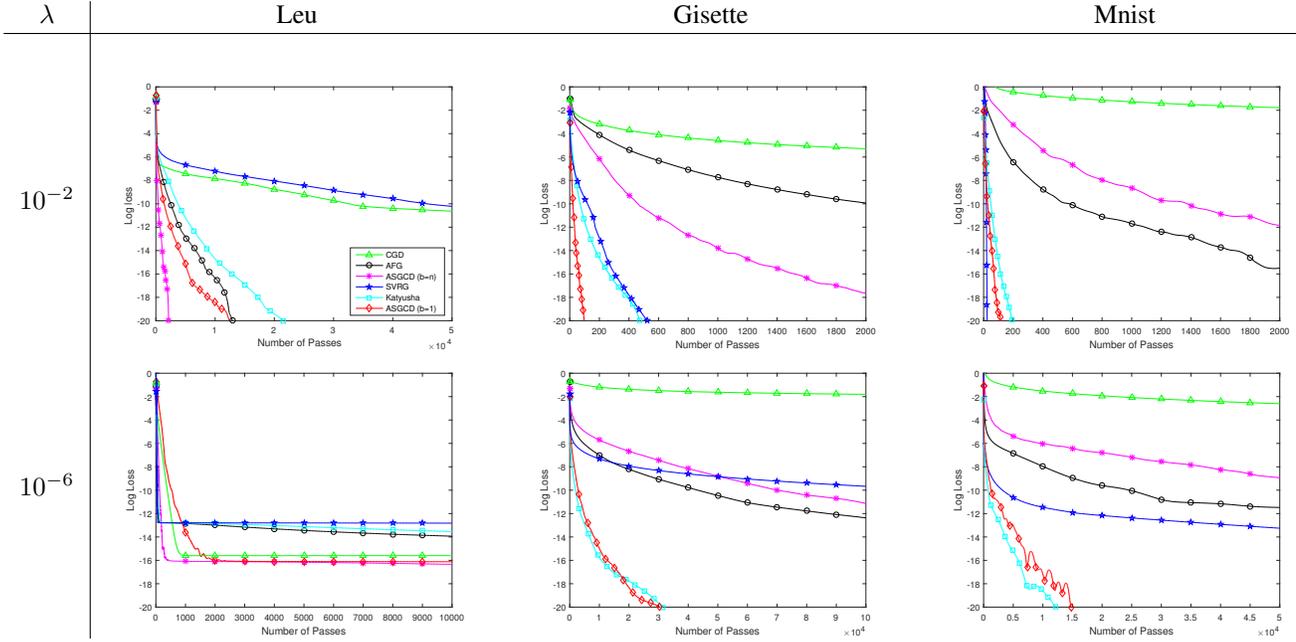

Fig. 1. Comparing AGCD ($b = 1$) and ASGCD ($b = n$) with CGD, SVRG, AFG and Katyusha on Lasso

TABLE II
FACTOR RATES OF FOR THE 6 CASES

| $\lambda$ | LEU | GISETTE | MNIST |
|---|---|---|---|
| $10^{-2}$ | $(0.85, 1.33)$ | $(0.88, 0.74)$ | $(5.85, 3.02)$ |
| $10^{-6}$ | $(1.45, 2.27)$ | $(3.51, 2.94)$ | $(5.84, 3.02)$ |

The performance of the 6 algorithms are plotted in Fig. 1. We use Log loss $\log(F(x_k) - F(x^*))$ in the $y$-axis. $x$-axis denotes the number that the algorithm access the data matrix $A$. For example, ASGCD with $b = n$ access $A$ once in each iteration, while ASGCD with $b = 1$ access $A$ twice in an entire outer iteration. For each case (Dataset, $\lambda$), we compute the rate $(r_1, r_2) = \left( \frac{\sqrt{CL_1}\|x^*\|_1}{\sqrt{L_2}\|x^*\|_2}, \frac{\sqrt{CT_1}\|x^*\|_1}{\sqrt{T_2}\|x^*\|_2} \right)$ in Table II. First, because of acceleration effect, ASGCD($b = n$) are always better than the non-accelerated CGD algorithm; second, compare ASGCD($b = 1$) with Katyusha and ASGCD($b = n$) with AFG, it is found that for the cases (Leu, $10^{-2}$), (Leu, $10^{-6}$) and (Gisette, $10^{-2}$), ASGCD($b = 1$) dominates Katyusha [1, Alg.2] and ASGCD($b = n$) dominates AFG. While the theoretical analysis in Section III show that if $r_1$ is relatively small such as around 1, then ASGCD($b = 1$) will be better than [1, Alg.2]. For the other 3 cases, [1, Alg.2] and AFG are better. The consistency between Table II and Fig. 1 demonstrates the theoretical analysis.

# APPENDIX

*Proof.* By using the variational identity $\|g\|_1^2 = \inf_{\theta \in \triangle_d} \sum_{i=1}^d \frac{g_i^2}{\theta_i}$ in [4] and the definition of $J(g, \theta)$ in (3), it follows that (3) can be rewritten as

$$\tilde{h} = \arg\min_{g \in \mathbb{R}^d} \{\inf_{\theta \in \triangle_d} J(g, \theta)\}.$$

By the joint convexity of $J(g, \theta)$, we can find the minimizer $\tilde{h}$ by swapping the optimization order of $g$ and $\theta$, which is to say based on the definition of $\tilde{g}(\theta)$, $J(\theta)$ and $\tilde{\theta}$,

$$\tilde{h} = \tilde{g}(\tilde{\theta}).$$

Therefore, the minimization problem of $\tilde{h}$ in (3) can be equivalently transformed to the problem (7). Meanwhile, it is observed that $J(g, \theta)$ is coordinate separable, *i.e.*,

$$J(g, \theta) = \sum_{i=1}^d J_i(g_i, \theta_i), \text{ where } J_i(g_i, \theta_i) \stackrel{\text{def}}{=} \nabla_i f(x) g_i + \frac{1}{2\eta} \frac{g_i^2}{\theta_i} + \lambda |x_i + g_i|.$$

By the definition of $\tilde{g}(\theta)$ in (6), $\tilde{g}(\theta)$ is also coordinate separable, *i.e.* for all $i \in [d]$,

$$\tilde{g}_i(\theta) = \tilde{g}_i(\theta_i) \stackrel{\text{def}}{=} \arg\min_{g_i \in \mathbb{R}} \left\{ \nabla_i f(x) g_i + \frac{1}{2\eta} \frac{g_i^2}{\theta_i} + \lambda |x_i + g_i| \right\}.$$

By using the iterative soft thresholding (IST) operator [6], for all $i \in [d]$,

$$\tilde{g}_i(\theta_i) = \text{sign}(x_i - \theta_i \eta \nabla_i f(x)) \cdot \max\{0, |x_i - \theta_i \eta \nabla_i f(x)| - \theta_i \eta \lambda\} - x_i.$$

Then it implies that $J(\theta)$ is also coordinate separable, *i.e.*,

$$J(\theta) = \sum_{i=1}^d J_i(\theta_i), \quad \text{where} \quad J_i(\theta_i) \stackrel{\text{def}}{=} J(\tilde{g}_i(\theta_i), \theta_i).$$

□

*Proof of Lemma 2.* For all $i \in [d]$ and $0 \leq \theta_i \leq 1$, by substituting (8) into (9), we get the expression of $J_i(\theta_i)$. Taking the derivative of $J_i(\theta_i)$ and set $\theta_i = 0, 1$ respectively, then we get the expressions of $J_i'(0), J_i'(1)$ as follows.

For all $i \in [d]$ and $0 \leq \theta_i \leq 1$, denote

$$\nu_i \stackrel{\text{def}}{=} -\frac{(\max\{|\nabla_i f(x)| - \lambda, 0\})^2 \eta}{2}, \qquad \chi_i(\theta_i) \stackrel{\text{def}}{=} -\frac{(\text{sign}(x_i - \theta_i \eta \nabla_i f(x))\lambda + \nabla_i f(x))^2 \eta}{2}, \tag{A.1}$$

then the derivate $J_i'(\theta_i)$ at $\theta_i = 0, 1$ are

$$J_i'(0) = \begin{cases} \nu_i, & x_i = 0, \\ \chi_i(0), & x_i \neq 0 \end{cases}, \quad J_i'(1) = \begin{cases} -\frac{x_i^2}{2\eta}, & |x_i - \eta \nabla_i f(x)| \leq \eta \lambda \\ \chi_i(1), & |x_i - \eta \nabla_i f(x)| > \eta \lambda \end{cases}. \tag{A.2}$$

For all $i \in [d]$, according to the value of $x_i$ and $\nabla_i f(x)$, by classified discussion, we can show that $J_i'(\theta_i)$ belongs to one of the 4 cases in Lemma 2. Firstly, we denote

$$O \stackrel{\text{def}}{=} \{i | 0 \in \nabla_i f(x) + \lambda \partial |x_i|\}, \tag{A.3}$$

$$\begin{aligned} U \stackrel{\text{def}}{=} & \{i \in [d] : x_i \geq 0, \nabla_i f(x) < -\lambda\} \\ & \cup \{i \in [d] | x_i \leq 0, \nabla_i f(x) > \lambda\} \\ & \cup \{i \in [d] | x_i > 0, \nabla_i f(x) > -\lambda, r_1 \geq 1\} \\ & \cup \{i \in [d] | x_i < 0, \nabla_i f(x) < \lambda, r_1 \geq 1\}, \end{aligned} \tag{A.4}$$

$$\begin{aligned} V \stackrel{\text{def}}{=} & \{i \in [d] | x_i > 0, -\lambda < \nabla_i f(x) \leq \lambda, r_1 < 1\} \\ & \cup \{i \in [d] | x_i < 0, -\lambda \leq \nabla_i f(x) < \lambda, r_1 < 1\} \\ & \cup \{i \in [d] | x_i > 0, \nabla_i f(x) > \lambda, r_2 \geq 1\} \\ & \cup \{i \in [d] | x_i < 0, \nabla_i f(x) < -\lambda, r_2 \geq 1\}, \end{aligned} \tag{A.5}$$

$$\begin{aligned} W \stackrel{\text{def}}{=} & \{i \in [d] | x_i > 0, \nabla_i f(x) > \lambda, r_2 < 1\} \\ & \cup \{i \in [d] | x_i < 0, \nabla_i f(x) < -\lambda, r_2 < 1\}. \end{aligned} \tag{A.6}$$

Then we can summarize the results as follows

- If $i \in O$, then $J'_i(\theta_i)$ belongs to the (case a) in Lemma 2.
- If $i \in U$, then $J'_i(\theta_i)$ belongs to the (case b) in Lemma 2.
- If $i \in V$, then $J'_i(\theta_i)$ belongs to the (case c) in Lemma 2.
- If $i \in W$, then $J'_i(\theta_i)$ belongs to the (case d) in Lemma 2.

$\square$

*Proof of Corollary 1.* Corollary 1 can be obtained by simply summarized the 4 cases in Lemma 2. $\square$

*Proof of Proposition 1.* Firstly, by checking $i \in O, U, V$ or $W$ orderly and using the expression of $J'_i(0)$ and $J'_i(1)$, it follows that for $i \in O \cup U$, $J'_i(0) = J'_i(1)$; for $i \in V$, $(-J'_i(0)) - (-J'_i(1)) \leq 2\lambda$; for $i \in W$, $(-J'_i(0)) - (-J'_i(1)) = 2\lambda$.

Therefore we have
$$\max_{i \in [d]}\{-J'_i(0)\} - \max_{k \in [d]}\{-J'_k(0)\} \leq \max_{i \in [d]}\{(-J'_i(0)) - (-J'_i(0))\} \leq 2\lambda.$$

and by Corollary 1,
$$0 \leq \max_{i \in [d]}\{-J'_i(0)\} - \max_{k \in [d]}\{-J'_k(0)\},$$

Proposition 1 is proved. $\square$

The simplex constraint is a linear constraint. By the property of KKT condition, if $(\tilde{\gamma}, \tilde{\theta}, \tilde{\zeta})$ is a stationary point of the problem (10), then $\tilde{\theta}$ is an optimal solution of (7). For $\hat{\theta}$, one can divide $[d]$ into two disjoint parts $S$ and $T$, where
$$S = \{i : \hat{\theta}_i > 0\} \text{ and } T = \{i : \hat{\theta}_i = 0\}.$$

Then $\forall i \in S$, by the complementary slackness $\hat{\zeta}_i \hat{\theta}_i = 0$, one has $\hat{\zeta}_i = 0$ and $\hat{\gamma} = -J'_i(\hat{\theta}_i) \geq 0$; $\forall j \in T$, similarly, one has $\hat{\zeta}_j \geq 0$ and $\hat{\gamma} \geq -J'_j(\hat{\theta}_j) \geq 0$. Thus the KKT condition can be equivalently written as (12).

*Proof of Lemma 4.* By Lemma 2, it follows that $J'_t(0) \leq 0$. Combing with the condition $J'_s(0) < J'_t(0)$, $J'_s(0) < 0$ and thus $J'_s(\theta_s)$ belongs to (case b), (case c) or (case d). Denote
$$r_s = \begin{cases} 1, & J'_s(\theta_s) \text{ belongs to (case b)}; \\ r_{s1}, & J'_s(\theta_s) \text{ belongs to (case c) or (case d)}, \end{cases} \tag{A.7}$$

where by the definition of $r_{i1}$ in Lemma 2, $r_{s1} = \frac{|x_s|}{\sqrt{-2\eta J'_s(0)}}$.

Assume by contradiction that $\tilde{\theta}_s = 0$ yet $\tilde{\theta}_t > 0$. Let $\hat{\theta}$ be a vector of which the elements are equal to the elements of $\tilde{\theta}$ except that
$$\hat{\theta}_s = \min\{\tilde{\theta}_t, r_s\}; \tag{A.8}$$
$$\hat{\theta}_t = \max\{0, \tilde{\theta}_s - r_s\}. \tag{A.9}$$

By the definition of $\hat{\theta}_s, \hat{\theta}_t$ in (A.8) and (A.9), it follows that
$$\forall \theta_s \in [\tilde{\theta}_s, \hat{\theta}_s], \quad J'_s(\theta_s) = J'_s(0) \tag{A.10}$$
$$\forall \theta_t \in [\hat{\theta}_t, \tilde{\theta}_t], \quad J'_t(\theta_t) \geq J'_t(0). \tag{A.11}$$

Then
$$\begin{aligned}
J(\tilde{\theta}) - J(\hat{\theta}) &= J_s(0) + \int_0^{\tilde{\theta}_s} J'_s(\theta_s)d\theta_s + J_s(0) + \int_0^{\tilde{\theta}_t} J'_t(\theta_t)d\theta_t \\
&\quad - J_s(0) - \int_0^{\hat{\theta}_s} J'_s(\theta_s)d\theta_s - J_t(0) - \int_0^{\hat{\theta}_t} J'_t(\theta_t)d\theta_t \\
&= \int_{\hat{\theta}_s}^{\tilde{\theta}_s} J'_s(\theta)d\theta + \int_{\hat{\theta}_t}^{\tilde{\theta}_t} J'_t(\theta) \\
&\geq \int_{\hat{\theta}_s}^{\tilde{\theta}_s} J'_s(0)d\theta_s + \int_{\hat{\theta}_t}^{\tilde{\theta}_t} J'_t(0)d\theta_t \\
&\geq J'_s(0)(\tilde{\theta}_s - \hat{\theta}_s) + J'_t(0)(\tilde{\theta}_t - \hat{\theta}_t)
\end{aligned}$$

Then the expressions of $\hat{\theta}_s, \hat{\theta}_t$ in (A.8) and (A.9),
$$J(\tilde{\theta}) - J(\hat{\theta}) = \begin{cases} (J'_t(0) - J'_s(0)) \cdot \tilde{\theta}_t, & \tilde{\theta}_t < r_s \\ (J'_t(0) - J'_s(0)) \cdot r_s, & \tilde{\theta}_t \geq r_s \end{cases}. \tag{A.12}$$

By the assumption $J'_s(0) < J'_t(0)$, $J(\tilde{\theta}) - J(\hat{\theta}) > 0$, which contradicts the fact that $\tilde{\theta}$ is the optimal solution. □

*Proof of Lemma 5.* By the KKT condition (12) in Lemma 3, it follows that for all $i \in S$, $-J'_i(\tilde{\theta}_i) \geq \max_{j \in T} -J'_j(0)$; meanwhile by Corollary 1, for all $i \in [d]$, $-J'_i(\tilde{\theta}_i)$ is a non-increasing function. Therefore combing the KKT condition (12), it follows that

$$\forall i \in S, \qquad -J'_i(0) \geq -J'_i(\tilde{\theta}_i) \geq \max_{j \in T}\{-J'_j(0)\} \geq \max_{j \in T}\{-J'_j(1)\} \tag{A.13}$$

In addition, by the KKT condition (12), for all $i_1 \in S, i_2 \in S$, $-J'_{i_1}(\tilde{\theta}_{i_1}) = -J'_{i_2}(\tilde{\theta}_{i_2})$. Because by Corollary 1, for all $i \in [d]$, $-J'_i(\tilde{\theta}_i)$ is a non-increasing function, therefore

$$\forall i_1 \in S, i_2 \in S, \qquad -J'_{i_1}(0) \geq -J'_{i_1}(\tilde{\theta}_{i_1}) = -J'_{i_2}(\tilde{\theta}_{i_2}) \geq -J'_{i_2}(1).$$

Therefore it follows that

$$\forall i \in S, \qquad -J'_i(0) \geq \max_{j \in S} -J'_j(1). \tag{A.14}$$

By combing (A.13) and (A.14), we get

$$-J'_i(0) \geq \max_{j \in [d]} -J'_j(1). \tag{A.15}$$

□

*Proof.* To prove Theorem 1, by Lemma 1, we only need to show $\tilde{\theta}$ in Alg. 1 is the optimal solution of the problem (7). By Lemma 3, to prove the optimality of $\tilde{\theta}$ in Alg. 1, we only need to show the $\tilde{\gamma}, \tilde{\theta}$ in Alg. 1 satisfy the KKT condition in Lemma 3. For convenience, we rewrite the KKT condition in this context,

$$\begin{cases} \sum_{i \in S} \tilde{\theta}_i = 1; & \text{(A.16a)} \\ \text{for all } j \in T, \quad \tilde{\theta}_j = 0; & \text{(A.16b)} \\ \text{for all } i \in S, \quad \tilde{\gamma} = -J'_i(\tilde{\theta}_i) \geq \max_{j \in T} -J'_j(0), & \text{(A.16c)} \end{cases}$$

where as in Lemma 3, $S = \{i \in [d] | \tilde{\theta}_i > 0\}, T = \{i \in [d] | \tilde{\theta}_i = 0\}$. The main difficulty in the proof of Theorem 1 comes from the fact that by Lemma 2, for all $i \in [d]$ and $0 \leq \theta_i \leq 1$, the expression of $J'_i(\theta_i)$ has 4 different cases. Here we give Lemma 6 to show an equivalence relation between the expression of $J'_i(\theta_i)$ and the relation of $J'_i(0)$ and $J'_i(1)$.

**Lemma 6.** *For all $i \in [d]$ and $0 \leq \theta_i \leq 1$, $J'_i(\theta_i)$ belongs to the (case a) or (case b) 2 if and only if $J'_i(0) = J'_i(1)$; $J'_i(\theta_i)$ belongs to the (case c) or (case d) in Lemma 2 if and only if $J'_i(0) \neq J'_i(1)$.*

*Proof of Lemma 6.* For all $i \in [d]$ and $0 \leq \theta_i 1$, by observing the (case a), (case b), (case c) and (case d) of $J'_i(\theta_i)$ in Lemma 2, it follows that $J'_i(\theta_i)$ belongs to the (case a) or (case b) if and only if it is a constant function, which implies $J'_i(0) = J'_i(1)$. $J'_i(\theta_i)$ belongs to the (case c) or (case d) if and only if it is a piecewise function, which implies $J'_i(0) \neq J'_i(1)$. □

By Lemma 6, the condition $J'_{i_j}(0) = J'_{i_j}(1)$ in the step 3 of Alg. 1 is used to identify which case $J'_{i_j}(\theta_{i_j})$ belongs to. Lemma 7 introduces an implied result of the conditions in the step 3 of Alg. 1. (In the following lemmas, we assume that $r_{i1}, r_{i2}(i \in [d])$ have been defined in Lemma 2 and $i_m, v_m, Q, v, \rho, \tilde{\gamma}, R_j(j \in [|Q|+1])$ have been defined in Alg. 1.)

**Lemma 7.** *For all $j \in [\rho - 1]$, it follows that all the following conditions*

$$\begin{cases} J'_{i_j}(0) \neq J'_{i_j}(1) & \text{(A.17a)} \\ \sum_{l \in R_j} |x_l| < \sqrt{2\eta v_j} & \text{(A.17b)} \\ j < |Q| + 1 & \text{(A.17c)} \end{cases}$$

*must be satisfied.*

*Proof of Lemma 7.* By the step 3 in Alg. 1, $\rho$ is the minimal index that can satisfies one of the following 3 condition

$$\begin{cases} v_\rho = -J'_{i_\rho}(1) \\ \sum_{l \in R_\rho} |x_l| \geq \sqrt{2\eta v_{i_\rho}} \\ \rho = |Q| + 1, \end{cases}$$

which implies that for all $j \in [\rho - 1]$, $j$ satisfies all the 3 conditions in Lemma 7. □

By Lemma 7, $j \in [\rho - 1]$ shares the 3 common properties in (A.17a)-(A.17c), which is important for the proof of the subsequent lemmas about $j \in [\rho - 1]$. In Lemma 8, we can find an useful inequalities.

**Lemma 8.** *For all $j \in [\rho - 1]$, $v_j \geq \tilde{\gamma} \geq v_\rho \geq v_m \geq -J'_{i_j}(1)$.*

*Proof of Lemma 8.* By the step 4 in Alg. 1, $\tilde{\gamma}$ has two possible values.

If $\tilde{\gamma} = (\sum_{k \in R_{\rho-1}} |x_k|)^2/(2\eta)$, it follows that
- In Lemma 7, let $j = \rho - 1$, we have $\sum_{k \in R_{\rho-1}} |x_k| < \sqrt{2\eta v_{\rho-1}}$. Then $v_{\rho-1} \geq \tilde{\gamma}$. Then by the definition of $v$, $v_1 \geq v_2 \geq \ldots \geq v_{\rho-1}$. Thus for all $j \in [\rho - 1]$, we have $v_j \geq \tilde{\gamma}$.
- In the step 4, when $\tilde{\gamma} = (\sum_{k \in R_{\rho-1}} |x_k|)^2/(2\eta)$, by the condition $\sum_{k \in R_{\rho-1}} |x_k| \geq \sqrt{2\eta v_\rho}$, we have $\tilde{\gamma} \geq v_\rho$.

If $\tilde{\gamma} = v_\rho$, by the definition of $v$, $v_1 \geq v_2 \geq \ldots \geq v_{\rho-1}$. Then it follows that for all $j \in [\rho - 1]$, $v_j \geq \tilde{\gamma} \geq v_\rho$.

By the definition of $v$, $v_1 \geq v_2 \geq \ldots \geq v_m$. By the definition of $v_m$, $v_m = \max_{i \in [d]} -J'_i(1) \geq -J'_{i'_j}(1)$. Therefore for all $j \in [\rho - 1]$, $v_\rho \geq v_m \geq -J'_{i_j}(1)$.

Combining the above analysis, Lemma 8 is proved. □

Lemma 8 gives $\tilde{\gamma}$ both lower and upper bounds, which then further bounds the range of $\tilde{\theta}_l = \frac{|x_l|}{\sqrt{2\eta\tilde{\gamma}}}$.

Before continue, we show the relation between $R_{\rho-1}$ and $R_\rho \setminus \{i_\rho\}$.

**Lemma 9.** *If $\rho < |Q| + 1$, then $R_{\rho-1} = R_\rho \setminus \{i_\rho\}$; if $v_{|Q|+1} = -J'_{i_{|Q|+1}}(1)$, then $R_{|Q|} = R_{|Q|+1} \setminus \{i_{|Q|+1}\}$; if $v_{|Q|+1} \neq -J'_{i_{|Q|+1}}(1)$, then $R_{|Q|} = R_{|Q|+1}$.*

*Proof of Lemma 9.*
- If $\rho < |Q| + 1$, then by the step 2 in Alg. 1, $i_1, i_2, \ldots, i_{|Q|}$ are different coordinates. Thus $R_{\rho-1} = R_\rho \setminus \{i_\rho\}$.
- If $v_{|Q|+1} = -J'_{i_{|Q|+1}}(1)$, by the definition of $Q$ in the step 1 of Alg. 1, $i_{|Q|+1} \notin Q$. Therefore $R_{|Q|} = R_{|Q|+1} \setminus \{i_{|Q|+1}\}$.
- If $v_{|Q|+1} \neq -J'_{i_{|Q|+1}}(1)$, then by Lemma 6, $J'_{i_{|Q|+1}}(\theta_{i_{|Q|+1}})$ belongs to (case c) or (case d). It follows that $-J'_{i_{|Q|+1}}(0) > -J'_{i_{|Q|+1}}(1) = v_m$. Then by the definition of $|Q|$ in the step 1 of Alg. 1, $i_{|Q|+1} \in Q$. Therefore $R_{|Q|} = R_{|Q|+1}$. □

In Lemma 10, by Lemma 8, we show that $J'_l(\tilde{\theta}_l) = -\tilde{\gamma}$, which is a part of the KKT condition (A.16c).

**Lemma 10.** *For all $l \in R_\rho \setminus \{i_\rho\}$, by setting $\tilde{\theta}_l = \frac{|x_l|}{\sqrt{2\eta\tilde{\gamma}}}$ as in the step 5 of Alg. 1, it follows that $J'_l(\tilde{\theta}_l) = -\tilde{\gamma}$.*

*Proof of Lemma 10.* By Lemma 9, $R_\rho \setminus \{i_\rho\} \subset R_{\rho-1}$. Then by Lemma 7, for all $l \in R_\rho \setminus \{i_\rho\}$, $J'_l(0) \neq J_l(1)$. Then by Lemma 6, for all $l \in R_\rho \setminus \{i_\rho\}$ and $0 \leq \theta_l \leq 1$, $J'_l(\theta_l)$ belongs to (case c) or (case d) in Lemma 2.

Assume that for $i \in [d]$, $r_{i1} = \frac{|x_i|}{\sqrt{-2\eta J'_i(0)}}$ and $r_{i1} = \frac{|x_i|}{\sqrt{-2\eta J'_i(1)}}$ have been defined in Lemma 2. By Lemma 8, we have $-J'_{i_j}(0) = v_j \geq \tilde{\gamma} \geq -J'_{i_j}(1)$. Then $r_{i1} \leq \tilde{\theta}_l = \frac{|x_l|}{\sqrt{2\eta\tilde{\gamma}}} \leq r_{i2}$.

In addition, if $\tilde{\gamma} = (\sum_{k \in R_{\rho-1}} |x_k|)^2/(2\eta)$, when $\tilde{\theta}_l = \frac{|x_l|}{\sum_{k \in R_{\rho-1}} |x_k|} \leq 1$; if $\tilde{\gamma} = v_\rho$, by the step 4 in Alg. 1, the condition $\sum_{k \in R_{\rho-1}} |x_k| \leq \sqrt{2\eta v_\rho}$ holds. Then $\tilde{\theta}_l = \frac{|x_k|}{\sqrt{2\eta\tilde{\gamma}}} = \frac{|x_k|}{\sqrt{2\eta v_\rho}} \leq 1$.

Therefore $r_{i1} \leq \tilde{\theta}_l \leq \min\{r_{i2}, 1\}$. By the form of (case c) and (case d) in Lemma 2, we can find that for all $l \in R_\rho \setminus \{i_\rho\}$, $J'_l(\tilde{\theta}_l) = -\frac{|x_l|^2}{2\eta\theta_l^2} = -\tilde{\gamma}$. □

To proof the KKT condition (A.16c), besides Lemma 10 for the case $j \in [\rho - 1]$, we also need Lemma 11 for the case $j = \rho$.

**Lemma 11.** *By setting $\tilde{\theta}_{i_\rho} = 1 - \sum_{k \in R_\rho \setminus \{i_\rho\}} \tilde{\theta}_k$ as in the step 5 of Alg. 1, then it follows that $\tilde{\theta}_{i_\rho} \geq 0$. Meanwhile, if $\tilde{\theta}_{i_\rho} = 0$, then $\tilde{\gamma} = (\sum_{k \in R_\rho \setminus \{i_\rho\}} |x_k|)^2/(2\eta)$; if $\tilde{\theta}_{i_\rho} > 0$, then $J'_{i_\rho}(\tilde{\theta}_{i_\rho}) = -\tilde{\gamma} = v_\rho$.*

*Proof of Lemma 11.* By the step 4 in Alg. 1, $\tilde{\gamma}$ has two possible values.

If $\tilde{\gamma} = (\sum_{k \in R_{\rho-1}} |x_k|)^2/(2\eta)$, we can give the analyses by discussing the following 3 cases.
- if $\rho < |Q| + 1$, then by Lemma 7, $R_\rho \setminus \{i_\rho\} = R_{\rho-1}$. By the step 5 in Alg. 1 and $\tilde{\gamma} = (\sum_{k \in R_{\rho-1}} |x_k|)^2/(2\eta)$, for all $l \in R_\rho \setminus \{i_\rho\}$, $\tilde{\theta}_l = \frac{|x_l|}{\sqrt{2\eta\tilde{\gamma}}} = \frac{|x_l|}{\sum_{k \in R_{\rho-1}} |x_l|}$. Then $\sum_{l \in R_\rho \setminus \{i_\rho\}} \tilde{\theta}_l = \sum_{l \in R_{\rho-1}} \tilde{\theta}_l = 1$. Therefore $\tilde{\theta}_{i_\rho} = 1 - \sum_{k \in R_\rho \setminus \{i_\rho\}} = 0$.
- If $\rho = |Q| + 1$ and $J'_{i_\rho}(0) = J'_{i_\rho}(1)$, then by Lemma 7, it follows that $R_\rho \setminus \{i_\rho\} = R_{\rho-1}$. Therefore, by the same analysis in the case $\rho < |Q| + 1$, $\tilde{\theta}_{i_\rho} = 1 - \sum_{k \in R_\rho \setminus \{i_\rho\}} = 0$.
- If $\rho = |Q| + 1$ and $J'_{i_\rho}(0) \neq J'_{i_\rho}(1)$, then by Lemma 7, $R_{\rho-1} = R_\rho = (R_\rho \setminus \{i_\rho\}) \cup \{i_\rho\}$. Therefore $0 \leq \tilde{\theta}_{i_\rho} = 1 - \sum_{k \in R_\rho \setminus \{i_\rho\}} \tilde{\theta}_k = \frac{|x_{i_\rho}|}{\sum_{k \in R_{\rho-1}} |x_k|} \leq 1$. Meanwhile, by Lemma 6, $J'_{i_\rho}(\tilde{\theta}_{i_\rho})$ belongs to (case c) or (case d) in Lemma 2. Due to $\tilde{\gamma} = (\sum_{k \in R_{\rho-1}} |x_k|)^2/(2\eta)$, then the condition $\sum_{k \in R_{\rho-1} |x_k|} \geq \sqrt{2\eta v_\rho}$ holds. Thus $\tilde{\theta}_{i_\rho} \leq \frac{|x_{i_\rho}|}{\sqrt{2\eta v_\rho}} = \frac{|x_{i_\rho}|}{\sqrt{2\eta v_{|Q|+1}}} =$ $\frac{|x_{i_\rho}|}{\sqrt{-2\eta J'_{i_\rho}(1)}} = r_{i_\rho 2}$, where $r_{i_\rho 2}$ is defined in Lemma 2. Meanwhile, in Lemma 7, let $j = \rho - 1$, then $\sum_{k \in R_{\rho-1}} |x_k| < \sqrt{2\eta v_\rho}$. Thus $\tilde{\theta}_{i_\rho} \geq \frac{|x_{i_\rho}|}{\sqrt{2\eta v_{\rho-1}}} = \frac{|x_{i_\rho}|}{\sqrt{2\eta v_{|Q|}}}$. Due to $i_\rho \in Q$, we have $v_{i_\rho} \geq v_{|Q|}$. $\tilde{\theta}_{i_\rho} \geq \frac{|x_{i_\rho}|}{\sqrt{2\eta v_\rho}} = \frac{|x_{i_\rho}|}{\sqrt{-2\eta J_{i_\rho}}} = r_{i_\rho 1}$, where $r_{i_\rho 1}$ is defined in Lemma 2. Combing the above analyses, we have $r_{i_\rho 1} \leq \tilde{\theta}_{i_\rho} \leq \min\{1, r_{i_\rho 2}\}$. Therefore by the form of (case c) and (case d) in Lemma 2, $J'_{i_\rho}(\tilde{\theta}_{i_\rho}) = -\frac{x_{i_\rho}^2}{2\eta \tilde{\theta}_{i_\rho}} = -\tilde{\gamma}$.

If $\tilde{\gamma} = v_\rho$, according to the condition in the step 4 of Alg. 1, $\sum_{l \in R_{\rho-1}} |x_l| < \sqrt{2\eta v_\rho}$. Meanwhile by Lemma 9, we have $R_{i_\rho} \setminus \{i_\rho\} \subset R_{i_\rho-1}$. Then $\sum_{l \in R_{i_\rho} \setminus \{i_\rho\}} \tilde{\theta}_l = \sum_{l \in R_{i_\rho} \setminus \{i_\rho\}} \frac{|x_l|}{\sqrt{2\eta \tilde{\gamma}}} \leq \sum_{l \in R_{i_\rho-1}} \frac{|x_l|}{\sqrt{2\eta \tilde{\gamma}}} < 1$. Therefore $\tilde{\theta}_{i_\rho} = 1 - \sum_{l \in R_{i_\rho} \setminus \{i_\rho\}} \tilde{\theta}_l > 0$. we can give the analyses by discussing the following 3 cases.

- If $v_\rho = -J'_{i_\rho}(1)$, then by Lemma 6, $J_{i_\rho}(\theta_{i_\rho})$ belongs to (case a) or (case b). Therefore $J_{i_\rho}(\tilde{\theta}_{i_\rho}) = J_{i_\rho}(0) = -v_\rho = \tilde{\gamma}$.
- If $\sum_{l \in R_\rho} |x_l| \geq \sqrt{2\eta v_\rho}$ and $v_\rho \neq -J'_{i_\rho}(1)$, then by Lemma 6, $J'_{i_\rho}(\theta_{i_\rho})$ belongs to (case c) or (case d) in Lemma 2. Meanwhile for $l \in R_\rho \setminus \{i_\rho\}$, by $\tilde{\theta}_l = \frac{|x_l|}{\sqrt{2\eta \tilde{\gamma}}} = \frac{|x_l|}{\sqrt{2\eta v_\rho}}$, we have $\tilde{\theta}_{i_\rho} = 1 - \sum_{l \in R_\rho \setminus \{i_\rho\}} \tilde{\theta}_l = 1 - \sum_{l \in R_\rho \setminus \{i_\rho\}} \frac{|x_l|}{\sqrt{2\eta v_\rho}} \leq \frac{|x_l|}{\sqrt{2\eta v_\rho}} = r_{i_\rho 1}$, where $r_{i_\rho 1}$ is defined in Lemma 2. Therefore by the form of (case c) and (case d) in Lemma 2, it follows that $J'_{i_\rho}(\tilde{\theta}_{i_\rho}) = J'_{i_\rho}(0) = -v_\rho = -\tilde{\gamma}$.
- If $\rho = |Q| + 1$, $\sum_{l \in R_\rho} |x_l| < \sqrt{2\eta v_\rho}$ and $v_\rho \neq -J'_{i_\rho}(1)$, then by $v_\rho \neq -J'_{i_\rho}(1)$ and Lemma 6, $J'_{i_\rho}(\theta_{i_\rho})$ belongs to (case c) or (case d) in Lemma 2. By $\sum_{l \in R_\rho} |x_l| < \sqrt{2\eta v_\rho}$, we have that the $r_{i_\rho 2}$ in Lemma 2 satisfies

$$r_{i_\rho 2} = \frac{|x_\rho|}{\sqrt{-2\eta J'_{i_\rho}}}(1) = \frac{|x_\rho|}{\sqrt{2\eta v_{|Q|+1}}} = \frac{|x_\rho|}{\sqrt{2\eta v_\rho}} < \frac{\sum_{l \in R_\rho} |x_l|}{\sqrt{2\eta v_\rho}} < 1.$$

Therefore $J'_{i_\rho}(\theta_{i_\rho})$ belongs to the (case d) in Lemma 2. For all $l \in R_\rho \setminus \{i_\rho\}$, by $\tilde{\theta}_l = \frac{|x_l|}{\sqrt{2\eta \tilde{\gamma}}} = \frac{|x_l|}{\sqrt{2\eta v_\rho}}$ and $\sum_{l \in R_\rho} |x_l| < \sqrt{2\eta v_\rho}$, we have

$$\tilde{\theta}_{i_\rho} = 1 - \sum_{k \in R_\rho \setminus \{i_\rho\}} \tilde{\theta}_k \geq \frac{|x_\rho|}{\sqrt{2\eta v_\rho}} = r_{i_\rho 2}. \tag{A.19}$$

Therefore, by the form of the (case d) in Lemma 2, $J'_{i_\rho}(\tilde{\theta}_{i_\rho} = J'_{i_\rho}(1) = -v_\rho = -\tilde{\gamma}$.

Summarizing the above analyses, Lemma 11 is proved.

$\square$

Then based on Lemmas 8, 10, 11 and the definition of $v_m$ in the step 1 of Alg. 1, we can show that the $\tilde{\gamma}, \tilde{\theta}$ obtained in Alg. 1 satisfies the KKT conditions (A.16a)-(A.16c).

- If $\tilde{\theta}_{i_\rho} = 0$, by Lemma 11, , then $\tilde{\gamma} = \sum_{k \in R_\rho \setminus \{i_\rho\}}$. Let $S = R_\rho \setminus \{i_\rho\}$ and $T = [d] \setminus (R_\rho \setminus \{i_\rho\})$. Let $S = R_\rho$ and $T = [d] \setminus R_\rho$, for all $i \in [d]$, by checking the value of $\tilde{\theta}_i$, it is found that the KKT conditions (A.16a)-(A.16c) are satisfied.
- If $\tilde{\theta}_{i_\rho} > 0$, by Lemma 11, $J'_{i_\rho}(\tilde{\theta}_{i_\rho}) = -\tilde{\gamma} = v_\rho$. Let $S = R_\rho$ and $T = [d] \setminus R_\rho$, for all $i \in [d]$, by checking the value of $\tilde{\theta}_i$, it is found that the KKT conditions (A.16a)-(A.16c) are satisfied.

Therefore Theorem 1 is proved. $\square$

### A. Some necessary Lemmas and Definitions

For $1 < p < \infty$ and the $\ell_p$-norm $\|\cdot\|_p$, we denote its dual norm as $\|x\|_q = \max_{\|y\|_p \leq 1} x^T y = (\sum_{i=1}^d |x_i|^q)^{\frac{1}{q}}$, where $\frac{1}{p} + \frac{1}{q} = 1$. For $p = 1$, by the definition of dual norm, the dual norm of $\ell_1$-norm is $\ell_\infty$-norm. In Lemma 12 and 13, some classical results are described.

**Lemma 12.** *If $\forall x, y \in \mathbb{R}^d, 1 \leq p \leq \infty, \frac{1}{p} + \frac{1}{q} = 1$, and $\eta > 0$, then $|\langle x, y \rangle| \leq \|x\|_q \|y\|_p \leq \frac{1}{2\eta} \|x\|_q^2 + \frac{\eta}{2} \|y\|_p^2$.*

**Lemma 13.** *If $\forall x \in \mathbb{R}^d, 1 \leq p \leq \infty$, then $\|x\|_p \leq \|x\|_1 \leq n^{1-\frac{1}{p}} \|x\|_p$.*

For a continuous differentiable function $f(x)$, we give the following definitions.

**Definition 1.** *$f(x)$ is $L_p$-smooth $(1 \leq p \leq \infty)$ w.r.t. $\|\cdot\|_p$ if $\forall x, y \in \mathbb{R}^d$ and $\frac{1}{p} + \frac{1}{q} = 1$, $\|\nabla f(x) - \nabla f(y)\|_q \leq L_p \|x - y\|_p$.*

By Definition. 1, we have Lemma 14.

**Lemma 14.** *If $f(x)$ is $L_p$-smooth $(1 \leq p \leq \infty)$ w.r.t $\|\cdot\|_p$ and $\frac{1}{p} + \frac{1}{q} = 1$, then*

$$\frac{1}{2L_p} \|\nabla f(x) - \nabla f(y)\|_q^2 \leq f(y) - f(x) - \langle \nabla f(x), y - x \rangle \leq \frac{L_p}{2} \|x - y\|_p^2. \tag{A.20}$$

*Proof.* Firstly it is showed that $f(x)$ being $L_p$-smooth w.r.t. $\|\cdot\|_p$ implies that $\forall x, y \in \mathbb{R}^d$

$$f(y) \leq f(x) + \langle \nabla f(x), y-x \rangle + \frac{L_p}{2}\|y-x\|_p^2.$$

Consider the function $g(\tau) = f(x+\tau(y-x))$ with $\tau \in \mathbb{R}$. Then

$$\begin{aligned} f(y) - f(x) - \langle \nabla f(x), y-x \rangle &= g(1) - g(0) - \langle \nabla f(x), y-x \rangle \\ &= \int_0^1 \frac{dg}{d\tau}(\tau) - \langle \nabla f(x), y-x \rangle \, d\tau \\ &= \int_0^1 \langle \nabla f(x+\tau(y-x)), y-x \rangle - \langle \nabla f(x), y-x \rangle \, d\tau \\ &= \int_0^1 \langle \nabla f(x+\tau(y-x)) - \nabla f(x), y-x \rangle \, d\tau \\ &\leq \int_0^1 \|\nabla f(x+\tau(y-x)) - \nabla f(x)\|_q \|y-x\|_p \, d\tau \\ &\leq \int_0^1 L_p \tau \|y-x\|_p^2 \, d\tau \\ &= \frac{L_p}{2}\tau^2 \|y-x\|_p^2 \Big|_0^1 \\ &= \frac{L_p}{2}\|y-x\|_p^2. \end{aligned}$$

To subsequently show $\frac{1}{2L_p}\|\nabla f(x) - \nabla f(y)\|_q^2 \leq f(y) - f(x) - \langle \nabla f(x), y-x \rangle$, fix $x \in \mathbb{R}^d$ and consider the function

$$\phi(y) = f(y) - \langle \nabla f(x), y \rangle,$$

which is convex on $\mathbb{R}^d$ and also has an $L_p$-Lipschitz continuous gradient w.r.t. $\|\cdot\|_p$, as

$$\begin{aligned} \|\phi'(y) - \phi'(x)\|_q &= \|(\nabla f(y) - \nabla f(x)) - (\nabla f(x) - \nabla f(x))\|_q \\ &= \|\nabla f(y) - \nabla f(x)\|_q \\ &\leq L_p \|y - x\|_p. \end{aligned}$$

As the minimizer of $\phi$ is $x$ (i.e., $\phi'(x) = 0$), for any $y \in \mathbb{R}^d$ we have

$$\begin{aligned} \phi(x) = \min_v \phi(v) &\leq \min_v \phi(y) + \langle \phi'(y), v-y \rangle + \frac{L_p}{2}\|v-y\|_p^2 \\ &= \phi(y) - \sup_v (\langle -\phi'(y), v-y \rangle - \frac{L_p}{2}\|v-y\|_p^2) \\ &= \phi(y) - \frac{1}{2L_p}\|\phi'(y)\|_q^2. \end{aligned}$$

Substituting in the definition of $\phi$, we have

$$f(x) - \langle \nabla f(x), x \rangle \leq f(y) - \langle \nabla f(x), y \rangle - \frac{1}{2L_p}\|\nabla f(y) - \nabla f(x)\|_q^2$$

$$\Longleftrightarrow \quad \frac{1}{2L_p}\|\nabla f(y) - \nabla f(x)\|_q^2 \leq f(y) - f(x) - \langle \nabla f(x), y-x \rangle.$$

□

**Definition 2.** *$f(x)$ is $\sigma_p$-strongly convex ($1 \leq p \leq \infty$) w.r.t. $\|\cdot\|_p$ if $\forall x, y \in \mathbb{R}^d$ and $\frac{1}{p} + \frac{1}{q} = 1$, $f(y) - f(x) - \langle \nabla f(x), y-x \rangle \geq \frac{\sigma_p}{2}\|x-y\|_p^2$.*

Taking $\frac{1}{2}\|x\|_p^2$ ($1 < p \leq 2$) as an example. It is known that $\frac{1}{2}\|x\|_p^2$ is $(p-1)$-strongly convex w.r.t. $\|\cdot\|_p$ [5]. Based on $\frac{1}{2}\|x\|_p^2$ ($1 < p \leq 2$), one can define $p$-Bregman divergence

$$B_p(y, x) = \frac{1}{2}\|y\|_p^2 - \frac{1}{2}\|x\|_p^2 - \langle \nabla \frac{1}{2}\|x\|_p^2, y-x \rangle. \tag{A.21}$$

**Lemma 15** ([5, 2]). *For $x, y \in \mathbb{R}^d, 1 < p \leq 2$, $B_p(y,x) = \frac{1}{2}\|y\|_p^2 - \frac{1}{2}\|x\|_p^2 - \langle \nabla \frac{1}{2}\|x\|_p^2, y-x \rangle$ satisfies the 3 properties.*
- $B_p(y, x) \geq \frac{p-1}{2}\|y-x\|_p^2$;
- $B_p(y, x) = 0$ if and only if $y = x$;
- $B_p(x, y) + B_p(y, z) = B_p(x, z) + \langle \frac{1}{2}\nabla\|z\|_p^2 - \frac{1}{2}\nabla\|y\|_p^2, x-y \rangle$.

## B. Proof of Theorem 2

Theorem 2 is proved by following the steps of the proof in [1]. Firstly, in Section B1, ASGCD is analyzed for a fixed $k$-th iteration. In the one-iteration analysis, $y_k, z_k$ and $x_{k+1}$ is assumed to be fixed and thus the selection of the mini batch $\mathcal{B}$ in the $k$-th iteration is the only source of randomness. For simplicity, let $\tilde{x} = \tilde{x}^s, \tau_{1,s} = \tau_1, \alpha = \alpha_s$ where $s = \lfloor \frac{k}{m} \rfloor$ is the epoch corresponding to $k$, let $\beta(b) = \frac{n-b}{b(n-1)}$ and denote $\sigma_{k+1}^2 = \|\nabla f(x_{k+1}) - \tilde{\nabla}_{k+1}\|_\infty^2$. Then $\mathbb{E}[\sigma_{k+1}^2]$ is the variance measured by $\|\cdot\|_\infty$ of the gradient estimator $\tilde{\nabla}_{k+1}$ in this iteration. Secondly, in Section B2, Theorem 2 is proven by combing the one-iteration analysis in Section B1 into the outer-iteration analysis in Section B2.

There are 3 differences from the analysis in [1]. Firstly, the analysis is used for the specific ASGCD algorithm that combines $\ell_1$-greedy and pComiD and thus the value of the parameter $\alpha_s$ is different from the value in [1]. Secondly, the analysis is given for mini batch selection based on $\|\cdot\|_\infty$ rather than one sample selection based on $\|\cdot\|_2$. Thirdly, we give the analysis for $\ell_1$-regularized problems (1) with a different way to represent result.

*1) One-iteration analysis:*

**Lemma 16** (SOTOPO). *If*

$$y_{k+1} = SOTOPO(\tilde{\nabla}_{k+1}, x_{k+1}, \lambda, \eta), \quad \text{and}$$
$$Prog(x_{k+1}) = -\min_{y \in \mathbb{R}^d} \left\{ \frac{(1+2\beta(b))L_1}{2} \|y - x_{k+1}\|_1^2 + \langle \tilde{\nabla}_{k+1}, y - x_{k+1} \rangle + \lambda\|y\|_1 - \lambda\|x_{k+1}\|_1 \right\} \geq 0,$$

*it follows that if $b < n$ (where the expectation is only over the randomness of $\tilde{\nabla}_{k+1}$), then*

$$F(x_{k+1}) - \mathbb{E}[F(y_{k+1})] \geq \mathbb{E}[Prog(x_{k+1})] - \frac{1}{4\beta(b)L_1}\mathbb{E}[\sigma_{k+1}^2]; \tag{A.22}$$

*if $b = n$ (no randomness exists), then*

$$F(x_{k+1}) - \mathbb{E}[F(y_{k+1})] \geq \mathbb{E}[Prog(x_{k+1})]. \tag{A.23}$$

*Proof.* If $b < n$, by setting $L = L_1$, it follows that

$$\begin{aligned}
\text{Prog}(x_{k+1}) &= -\min_y \left\{ \frac{(1+2\beta(b))L_1}{2} \|y - x_{k+1}\|_1^2 + \langle \tilde{\nabla}_{k+1}, y - x_{k+1} \rangle + \lambda\|y\|_1 - \lambda\|x_{k+1}\|_1 \right\} \\
&\stackrel{①}{=} -\frac{(1+2\beta(b))L_1}{2} \|y_{k+1} - x_{k+1}\|_1^2 + \langle \tilde{\nabla}_{k+1}, y_{k+1} - x_{k+1} \rangle + \lambda\|y_{k+1}\|_1 - \lambda\|x_{k+1}\|_1 \\
&= -\frac{L_1}{2} \|y_{k+1} - x_{k+1}\|_1^2 + \langle \nabla f(x_{k+1}), y_{k+1} - x_{k+1} \rangle + \lambda\|y_{k+1}\|_1 - \lambda\|x_{k+1}\|_1 \\
&\quad + \left( \langle \nabla f(x_{k+1}) - \tilde{\nabla}_{k+1}, y_{k+1} - x_{k+1} \rangle - \beta(b)L\|y_{k+1} - x_{k+1}\|_1^2 \right) \\
&\stackrel{②}{\leq} -(f(y_{k+1}) - f(x_{k+1}) + \lambda\|y_{k+1}\|_1 - \lambda\|x_{k+1}\|_1) + \frac{1}{4\beta(b)L_1} \|\nabla f(x_{k+1}) - \tilde{\nabla}_{k+1}\|_\infty^2,
\end{aligned}$$

where ① is by Theorem 1, ② is by the smoothness assumption (2), Lemma 12 and 14. Taking expectation on both sides, (A.22) is obtained.

If $b = n$, then $\beta(b) = 0$. By using a similar analysis as the case $b < n$, (A.23) is obtained. $\square$

**Lemma 17.** *(variance upper bound). If $b < n$, then*

$$\mathbb{E}[\|\tilde{\nabla}_{k+1} - \nabla f(x_{k+1})\|_\infty^2] \leq 2\beta(b)L_1(f(\tilde{x}) - f(x_{k+1}) - \langle \nabla f(x_{k+1}), \tilde{x} - x_{k+1} \rangle). \tag{A.24}$$

*Proof.* Before the proof, it should be noted that the variance upper bound measured by $\|\cdot\|_2$ of mini-batch selection has been proved in [23]. The variance in our case is measured by $\|\cdot\|_\infty$. Because some properties of $\|\cdot\|_2$ such as $\mathbb{E}[\|x - \mathbb{E}[x]\|_2^2] = \mathbb{E}[\|x\|_2^2] - \|\mathbb{E}[x]\|_2^2$ and $\|\sum_i x_i\|_2^2 = \sum_{i,j} x_i^T x_j$ can't be generalized to $\|\cdot\|_\infty$ directly, the proof is slightly different from the proof in [23].

Let $\phi_j = (\nabla f_j(x_{k+1}) - \nabla f_j(\tilde{x})) - (\nabla f(x_{k+1}) - \nabla f(\tilde{x}))$ and $\phi_j^i = (\nabla_i f_j(x_{k+1}) - \nabla_i f_j(\tilde{x})) - (\nabla_i f(x_{k+1}) - \nabla_i f(\tilde{x}))$. Denote $i_{\max} = \arg\max_i |\nabla_i f(x_{k+1}) - \tilde{\nabla}_{k+1,i}|$. It follows that

$$\begin{aligned}
\mathbb{E}\left[\left\|\frac{1}{b}\sum_{j\in\mathcal{B}}\phi_j^{i_{\max}}\right\|_\infty^2\right] &= \frac{1}{b^2}\mathbb{E}\left[\sum_{j_1,j_2\in\mathcal{B}}\phi_{j_1}^{i_{\max}}\phi_{j_2}^{i_{\max}}\right] \\
&= \frac{1}{b^2}\mathbb{E}\left[\sum_{j_1\neq j_2\in\mathcal{B}}\phi_{j_1}^{i_{\max}}\phi_{j_2}^{i_{\max}}\right] + \frac{1}{b}\mathbb{E}\left[(\phi_j^{i_{\max}})^2\right] \\
&= \frac{b-1}{bn(n-1)}\sum_{j_1\neq j_2\in[n]}\phi_{j_1}^{i_{\max}}\phi_{j_2}^{i_{\max}} + \frac{1}{b}\mathbb{E}\left[(\phi_j^{i_{\max}})^2\right] \\
&= \frac{b-1}{bn(n-1)}\sum_{j_1,j_2\in[n]}\phi_{j_1}^{i_{\max}}\phi_{j_2}^{i_{\max}} - \frac{b-1}{b(n-1)}\mathbb{E}\left[(\phi_j^{i_{\max}})^2\right] + \frac{1}{b}\mathbb{E}\left[(\phi_j^{i_{\max}})^2\right] \\
&= \frac{b-1}{bn(n-1)}\sum_{j_1,j_2\in[n]}\phi_{j_1}^{i_{\max}}\phi_{j_2}^{i_{\max}} - \beta(b)\mathbb{E}\left[(\phi_j^{i_{\max}})^2\right] \\
&\stackrel{①}{=} \beta(b)\mathbb{E}\left[(\phi_j^{i_{\max}})^2\right] \\
&\stackrel{②}{\leq} \beta(b)\mathbb{E}\left[\|\phi_j\|_\infty^2\right],
\end{aligned} \tag{A.25}$$

where ① is using the fact $\sum_{j=1}^n \phi_j^{i_{\max}} = 0$. ② is by definition of $\|\cdot\|_\infty$. Denote $i_{j_{\max}} = \arg\max_i |\phi_j^i|$. Hence

$$\begin{aligned}
&\mathbb{E}\left[\left\|\nabla f(x_{k+1}) - \tilde{\nabla}_{k+1}\right\|_\infty^2\right] \\
&= \mathbb{E}\left[\left\|\frac{1}{b}\sum_{j\in\mathcal{B}}(\nabla f_j(x_{k+1}) - \nabla f_j(\tilde{x})) - (\nabla f(x_{k+1}) - \nabla f(\tilde{x}))\right\|_\infty^2\right] \\
&\stackrel{①}{=} \beta(b)\mathbb{E}\left[\|\nabla f_j(x_{k+1}) - \nabla f_j(\tilde{x})) - (\nabla f(x_{k+1}) - \nabla f(\tilde{x})\|_\infty^2\right] \\
&\stackrel{②}{=} \beta(b)\mathbb{E}\left[\left(\nabla_{i_{j_{\max}}} f_j(x_{k+1}) - \nabla_{i_{j_{\max}}} f_j(\tilde{x})) - \nabla_{i_{j_{\max}}} f(x_{k+1}) - \nabla_{i_{j_{\max}}} f(\tilde{x})\right)^2\right] \\
&= \beta(b)\mathbb{E}\left[\left(\nabla_{i_{j_{\max}}} f_j(x_{k+1}) - \nabla_{i_{j_{\max}}} f_j(\tilde{x})\right)^2 - \left(\nabla_{i_{j_{\max}}} f(x_{k+1}) - \nabla_{i_{j_{\max}}} f(\tilde{x})\right)^2\right] \\
&\leq \beta(b)\mathbb{E}\left[\left(\nabla_{i_{j_{\max}}} f_j(x_{k+1}) - \nabla_{i_{j_{\max}}} f_j(\tilde{x})\right)^2\right] \\
&\stackrel{③}{\leq} \beta(b)\mathbb{E}\left[\|\nabla f_j(x_{k+1}) - \nabla f_j(\tilde{x})\|_\infty^2\right] \\
&\leq 2\beta(b)L_1\mathbb{E}\left[f_j(\tilde{x}) - f_j(x_{k+1}) - \langle \nabla f_j(x_{k+1}), \tilde{x} - x_{k+1}\rangle\right] \\
&\stackrel{④}{=} 2\beta(b)L_1(f(\tilde{x}) - f(x_{k+1}) - \langle \nabla f(x_{k+1}), \tilde{x} - x_{k+1}\rangle),
\end{aligned}$$

where ① is by (A.25), ② is using the fact $\mathbb{E}[(x - \mathbb{E}[x])] = \mathbb{E}[x^2] - (\mathbb{E}[x])^2$, ④ is by the definition of $\|\cdot\|_\infty$, ③ is by Lemma 14. □

**Lemma 18** (pCOMID). *Fixing $\tilde{\nabla}_{k+1}$ and letting*

$$(z_{k+1}, \theta_{k+1}) = pCOMID(\tilde{\nabla}_{k+1}, \theta_k, q, \lambda, \alpha), \tag{A.26}$$

*it satisfies for all $u \in \mathbb{R}^d$,*

$$\alpha\langle\tilde{\nabla}_{k+1}, z_{k+1} - u\rangle + \alpha\lambda\|z_{k+1}\|_1 - \alpha\lambda\|u\|_1 \leq -B_p(z_{k+1}, z_k) + B_p(u, z_k) - B_p(u, z_{k+1}). \tag{A.27}$$

*Proof.* From [11], we known that pCOMID exactly solves the following mirror descent problem,

$$z_{k+1} = \arg\min_z\{\langle\tilde{\nabla}_{k+1}, z - z_k\rangle + \frac{1}{\alpha}B_p(z, z_k) + \lambda\|z\|_1\}. \tag{A.28}$$

By the optimality condition of $z_{k+1}$, it follows that

$$\nabla\frac{1}{2}\|z_{k+1}\|_p^2 - \nabla\frac{1}{2}\|z_k\|_p^2 + \alpha\tilde{\nabla}_{k+1} + \alpha g = 0.$$

Then the equality

$$\langle \nabla \frac{1}{2}\|z_{k+1}\|_p^2 - \nabla \frac{1}{2}\|z_k\|_p^2 + \alpha\tilde{\nabla}_{k+1} + \alpha g, z_{k+1} - u\rangle = 0$$

holds. In addition by Lemma 15, it follows that $\langle \nabla \frac{1}{2}\|z_{k+1}\|_p^2 - \nabla \frac{1}{2}\|z_k\|_p^2, z_{k+1} - u\rangle = B_p(z_{k+1}, z_k) - B_p(u, z_k) + B_p(u, z_{k+1})$. By the convexity of $\lambda\|z\|_1$, $\langle g, z_{k+1} - u\rangle \geq \lambda\|z_{k+1}\|_1 - \lambda\|u\|_1$. Therefore, we can write

$$\begin{aligned}
&\alpha\langle \tilde{\nabla}_{k+1}, z_{k+1} - u\rangle + \alpha\lambda\|z_{k+1}\|_1 - \alpha\lambda\|u\|_1 \\
=& -\langle \nabla \frac{1}{2}\|z_{k+1}\|_p^2 - \nabla \frac{1}{2}\|z_k\|_p^2, z_{k+1} - u\rangle - \langle \alpha g, z_{k+1} - u\rangle + \alpha\lambda\|z_{k=1}\|_1 - \alpha\lambda\|u\|_1 \\
\leq& -B_p(z_{k+1}, z_k) + B_p(u, z_k) - B_p(u, z_{k+1}).
\end{aligned}$$

$\square$

**Lemma 19** (Couping step 1). *If $x_{k+1} = \tau_1 z_k + \tau_2 \tilde{x} + (1 - \tau_1 - \tau_1) y_k$, where $\tau_1 \leq \frac{1}{(1+\frac{2(n-b)}{b(n-1)})\alpha L_1}$ and $\tau_2 = \frac{1}{2}$,*

$$\begin{aligned}
&\alpha\langle \nabla f(x_{k+1}), z_k - u\rangle - \alpha\lambda\|u\|_1 \\
\leq& \frac{\alpha}{\tau_1}\left(F(x_{k+1}) - \mathbb{E}[F(y_{k+1})]\right) + \tau_2 F(\tilde{x}) - \tau_2 \mathbb{E}[F(x_{k+1})] - \tau_2\langle \nabla f(x_{k+1}), \tilde{x} - x_{k+1}\rangle) \\
&+ B_p(u, z_k) - \mathbb{E}[B_p(u, z_{k+1})] + \frac{\alpha(1-\tau_1-\tau_2)}{\tau_1}\lambda\|y_k\|_2 - \frac{\alpha}{\tau_1}\lambda\|x_{k+1}\|_1.
\end{aligned}$$

*Proof.* It follows taht

$$\begin{aligned}
&\alpha\langle \tilde{\nabla}_{k+1}, z_k - u\rangle + \alpha\lambda\|z_{k+1}\|_1 - \alpha\lambda\|u\|_1 \\
=& \alpha\langle \tilde{\nabla}_{k+1}, z_k - z_{k+1}\rangle + \alpha\langle \tilde{\nabla}_{k+1}, z_{k+1} - u\rangle + \alpha\lambda\|z_{k+1}\|_1 - \alpha\lambda\|u\|_1 \\
\overset{\text{①}}{\leq}& \alpha\langle \tilde{\nabla}_{k+1}, z_k - z_{k+1}\rangle - B_p(z_{k+1}, z_k) + B_p(u, z_k) - B_p(u, z_{k+1}) \\
\overset{\text{②}}{\leq}& \alpha\langle \tilde{\nabla}_{k+1}, z_k - z_{k+1}\rangle - \frac{p-1}{2}\|z_{k+1} - z_k\|_p^2 + B_p(u, z_k) - B_p(u, z_{k+1}) \\
\overset{\text{③}}{\leq}& \alpha\langle \tilde{\nabla}_{k+1}, z_k - z_{k+1}\rangle - \frac{p-1}{2}d^{-(1-\frac{1}{p})}\|z_{k+1} - z_k\|_1^2 + B_p(u, z_k) - B_p(u, z_{k+1}), \\
\overset{\text{④}}{=}& \alpha\langle \tilde{\nabla}_{k+1}, z_k - z_{k+1}\rangle - \frac{1}{2C}\|z_{k+1} - z_k\|_1^2 + B_p(u, z_k) - B_p(u, z_{k+1}),
\end{aligned} \quad (A.29)$$

where ① is by Lemma 18, ② is by Lemma 15, ③ is by Lemma 13 and ④ is by the setting $C = \frac{d^{\frac{2\delta}{1+\delta}}}{\delta}$ and $p = 1 + \delta$ in Alg. III.

By defining $v \overset{def}{=} \tau_1 z_{k+1} + \tau_2 \tilde{x} + (1 - \tau_1 - \tau_2) y_k$, we have $x_{k+1} - v = \tau_1(z_k - z_{k+1})$ and therefore

$$\begin{aligned}
&\mathbb{E}[\alpha\langle \tilde{\nabla}_{k+1}, z_k - z_{k+1}\rangle - \frac{1}{2C}\|z_{k+1} - z_k\|_1^2] = \mathbb{E}[\frac{\alpha}{\tau_1}\langle \tilde{\nabla}_{k+1}, x_{k+1} - v\rangle - \frac{1}{2C\tau_1^2}\|x_{k+1} - v\|_1^2] \\
=& \mathbb{E}\left[\frac{\alpha}{\tau_1}\left(\langle \tilde{\nabla}_{k+1}, x_{k+1} - v\rangle - \frac{1}{2C\alpha\tau_1}\|x_{k+1} - v\|_1^2 - \lambda\|v\|_1 + \lambda\|x_{k+1}\|_1\right) + \frac{\alpha}{\tau_1}\left(\lambda\|v\|_1 - \lambda\|x_{k+1}\|_1\right)\right] \\
\overset{\text{①}}{=}& \mathbb{E}\left[\frac{\alpha}{\tau_1}\left(\langle \tilde{\nabla}_{k+1}, x_{k+1} - v\rangle - \frac{(1+2\beta(b))L}{2}\|x_{k+1} - v\|_1^2 - \lambda\|v\|_1 + \lambda\|x_{k+1}\|_1\right) + \frac{\alpha}{\tau_1}\left(\lambda\|v\|_1 - \lambda\|x_{k+1}\|_1\right)\right] \\
\overset{\text{②}}{\leq}& \mathbb{E}\left[\frac{\alpha}{\tau_1}(F(x_{k+1}) - F(y_{k+1}) + \frac{1}{4\beta(b)L_1}\sigma_{k+1}^2) + \frac{\alpha}{\tau_1}(\lambda\|v\|_1 - \lambda\|x_{k+1}\|_1)\right] \\
\overset{\text{③}}{\leq}& \mathbb{E}\Big[\frac{\alpha}{\tau_1}\left(F(x_{k+1}) - F(y_{k+1}) + \frac{1}{2}(f(\tilde{x}) - f(x_{k+1}) - \langle \nabla f(x_{k+1}, \tilde{x} - x_{k+1})\rangle)\right) \\
&+ \frac{\alpha}{\tau_1}(\tau_1\lambda\|z_{k+1}\|_1 + \tau_2\lambda\|\tilde{x}\|_1 + (1-\tau_1-\tau_2)\lambda\|y_k\|_1 - \lambda\|x_{k+1}\|_1)\Big],
\end{aligned} \quad (A.30)$$

where ① is by the setting $\alpha\tau_1 = \frac{1}{(1+2\beta(b))CL}$, ② is by Lemma 16, ③ is by Lemma 17 and the convexity $\|v\|_1 = \|\tau_1 z_{k+1} + \tau_2 \tilde{x} + (1-\tau_1-\tau_2)y_k\|_1 \leq \tau_1\|z_{k+1}\|_1 + \tau_2\|\tilde{x}\|_1 + (1-\tau_1-\tau_2)\|y_k\|_1$. Then, it is showed that $\mathbb{E}[\langle \tilde{\nabla}_{k+1}, z_k - u\rangle] = \langle \nabla f(x_{k+1}), z_k - u\rangle$ and $\tau_2 = \frac{1}{2}$. Combing (A.29) and (A.30), Lemma 19 is obtained. $\square$

**Lemma 20** (Coupling step 2). *Under the same choices of $\tau_1, \tau_2$ as in Lemma 19, one has*

$$0 \leq \frac{\alpha(1-\tau_1-\tau_2)}{\tau_1}(F(y_k) - F(x^*)) - \frac{\alpha}{\tau_1}(\mathbb{E}[F(y_{k+1})] - F(x^*)) + \frac{\alpha\tau_2}{\tau_1}(F(\tilde{x}) - F(x^*)) + B_p(x^*, z_k) - \mathbb{E}[B_p(x^*, z_{k+1})].$$

*Proof.* It follows that

$$\begin{aligned}
\alpha(f(x_{k+1}) - f(u)) &\overset{①}{\leq} \alpha\langle f(x_{k+1}), x_{k+1} - u\rangle \\
&= \alpha\langle f(x_{k+1}), x_{k+1} - z_k\rangle + \alpha\langle \nabla f(x_{k+1}), z_k - u\rangle \\
&\overset{②}{=} \frac{\alpha\tau_2}{\tau_1}\langle \nabla f(x_{k+1}), \tilde{x} - x_{k+1}\rangle + \frac{\alpha(1-\tau_1-\tau_2)}{\tau_1}\langle \nabla f(x_{k+1}), y_k - x_{k+1}\rangle + \alpha\langle \nabla f(x_{k+1}), z_k - u\rangle \\
&\overset{③}{\leq} \frac{\alpha\tau_2}{\tau_1}\langle \nabla f(x_{k+1}), \tilde{x} - x_{k+1}\rangle + \frac{\alpha(1-\tau_1-\tau_2)}{\tau_1}(f(y_k) - f(x_{k+1})) + \alpha\langle f(x_{k+1}), z_k - u\rangle,
\end{aligned} \quad (A.31)$$

where ① is by the convexity of $f(x)$, ② is by the convex combination $x_{k+1} = \tau_1 z_k + \tau_2 \tilde{x} + (1 - \tau_1 - \tau_2)y_k$, ③ is again by the convexity of $f(x)$. Applying Lemma 19 to (A.31), it follows that

$$\alpha(f(x_{k+1}) - F(u)) \leq \frac{\alpha(1-\tau_1-\tau_2)}{\tau_1}(F(y_k) - f(x_{k+1}))$$
$$+ \frac{\alpha}{\tau_1}(F(x_{k+1}) - \mathbb{E}[F(y_{k+1})] + \tau_2 F(\tilde{x}) - \tau_2 f(x_{k+1})) + B_p(u, z_k) - \mathbb{E}[B_p(u, z_{k+1})] - \frac{\alpha}{\tau_1}\lambda\|x_{k+1}\|_1,$$

which implies

$$\alpha(F(x_{k+1}) - F(u)) \leq \frac{\alpha(1-\tau_1-\tau_2)}{\tau_1}(F(y_k) - F(x_{k+1}))$$
$$+ \frac{\alpha}{\tau_1}(F(x_{k+1}) - \mathbb{E}[F(y_{k+1})] + \tau_2 F(\tilde{x}) - \tau_2 F(x_{k+1})) + B_p(u, z_k) - \mathbb{E}[B_p(u, z_{k+1})].$$

After arrangement and setting $u$ to some minimizer $x^*$, Lemma 20 is obtained. □

*2) Proof of Theorem 2:*

*Proof.* Assume the parameter $\tau_{1,s}$ and $\alpha_s$ satisfies the assumption $\tau_{1,s}\alpha_s = \frac{1}{(1+2\beta(b))CL_1}$ in Lemma 19. Let $D_k \overset{def}{=} F(y_k) - F(x^*)$ and $\tilde{D}^s \overset{def}{=} F(\tilde{x}^s) - F(x^*)$, Lemma 20 can be rewritten as

$$0 \leq \frac{\alpha_s(1-\tau_{1,s}-\tau_2)}{\tau_{1,s}}D_k - \frac{\alpha_s}{\tau_{1,s}}\mathbb{E}[D_{k+1}] + \frac{\alpha_s\tau_2}{\tau_{1,s}}\tilde{D}^s + B_p(x^*, z_k) - \mathbb{E}[B_p(x^*, z_{k+1})].$$

In the $s$-th epoch, summing up the above inequality for all the iterations $k = sm, sm+1, \ldots, sm+m-1$, it follows that

$$\mathbb{E}\left[\alpha_s \frac{1-\tau_{1,s}-\tau_2}{\tau_{1,s}}D_{(s+1)m} + \alpha_s \frac{\tau_{1,s}+\tau_2}{\tau_{1,s}}\sum_{l=1}^{m} D_{sm+l}\right]$$
$$\leq \alpha_s \frac{1-\tau_{1,s}-\tau_2}{\tau_{1,s}}D_{sm} + \alpha_s \frac{\tau_2}{\tau_{1,s}}m\tilde{D}^s + B_p(x^*, z_{sm}) - \mathbb{E}[B_p(x^*, z_{(s+1)m})]. \quad (A.32)$$

It should be noted that in (A.32), we fix all the randomness in the first $s-1$ epochs and take expectation on the current epoch $s$.

By the definition $\tilde{x}^s = \frac{1}{m}\sum_{l=1}^{m} y_{(s-1)m+l}$ in Alg. III and the convexity of $F(x)$, we have $m\tilde{D}^s \leq \sum_{l=1}^{m} D_{(s-1)m+l}$. Then for each $s \geq 1$, by (A.32), it follows that

$$\mathbb{E}\left[\frac{1}{\tau_{1,s}^2}D_{(s+1)m} + \frac{\tau_{1,s}+\tau_2}{\tau_{1,s}^2}\sum_{l=1}^{m-1} D_{sm+l}\right]$$
$$\leq \frac{1-\tau_{1,s}}{\tau_{1,s}^2}D_{sm} + \frac{\tau_2}{\tau_{1,s}^2}\sum_{l=1}^{m-1} D_{(s-1)m+l} + (1+2\beta(b))CLB_p(x^*, z_{sm}) - (1+2\beta(b))CL\mathbb{E}[B_p(x^*, z_{(s+1)m})].$$

For $s = 0$, (A.32) can be written as

$$\mathbb{E}\left[\frac{1}{\tau_{1,0}^2}D_m + \frac{\tau_{1,0}+\tau_2}{\tau_{1,0}^2}\sum_{l=1}^{m-1} D_l\right]$$
$$\leq \frac{1-\tau_{1,0}-\tau_2}{\tau_{1,0}^2}D_0 + \frac{\tau_2 m}{\tau_{1,0}^2}\tilde{D}^0 + (1+2\beta(b))CLB_p(x^*, z_0) - (1+2\beta(b))CL\mathbb{E}[B_p(x^*, z_m)]. \quad (A.33)$$

Choose $\tau_{1,s} = \frac{2}{s+4} \leq \frac{1}{2}$ which satisfies

$$\frac{1}{\tau_{1,s}^2} \geq \frac{1-\tau_{1,s+1}}{\tau_{1,s+1}^2} \quad \text{and} \quad \frac{\tau_{1,s}+\tau_2}{\tau_{1,s}^2} \geq \frac{\tau_2}{\tau_{1,s+1}^2}. \quad (A.34)$$

Then it follows that

$$\mathbb{E}\left[\frac{m}{\tau_{1,S-1}^2}\tilde{D}^S + (1+2\beta(b))\,CLB_p(x^*, z_{Sm})\right]$$

$$\overset{\text{\textcircled{1}}}{\leq} \mathbb{E}\left[\frac{1}{\tau_{1,S-1}^2}D_{Sm} + \frac{\tau_{1,S-1}}{\tau_{1,S-1}^2}\sum_{l=1}^{m-1}D_{(S-1)m+l} + (1+2\beta(b))\,CLB_p(x^*, z_{Sm})\right]$$

$$\overset{\text{\textcircled{2}}}{\leq} \mathbb{E}\left[\frac{1}{\tau_{1,S-1}^2}D_{Sm} + \frac{\tau_2}{\tau_{1,S}^2}\sum_{l=1}^{m-1}D_{(S-1)m+l} + (1+2\beta(b))\,CLB_p(x^*, z_{Sm})\right]$$

$$\overset{\text{\textcircled{3}}}{\leq} \frac{1-\tau_{1,0}-\tau_2}{\tau_{1,0}^2}D_0 + \frac{\tau_2 m}{\tau_{1,0}^2}\tilde{D}^0 + (1+2\beta(b))\,CLB_p(x^*, z_0)$$

$$\overset{\text{\textcircled{4}}}{=} \frac{\tau_2 m}{\tau_{1,0}^2}\tilde{D}^0 + (1+2\beta(b))\,CLB_p(x^*, z_0),$$

where \textcircled{1} is by $m\tilde{D}^s \leq \sum_{l=1}^{m}D_{(s-1)m+l}$, \textcircled{2} is by $\tau_2 \geq \tau_{1,S-1} \geq \tau_{1,S}$, \textcircled{3} uses (A.34) to telescope (A.32) and (A.33) for all $s = 0, 1, \ldots, S-1$ and \textcircled{4} is by $\tau_{1,0} = \tau_2 = \frac{1}{2}$.

$$\begin{aligned}\mathbb{E}[F(\tilde{x}^S) - F(x^*)] &= \mathbb{E}\tilde{D}^S \leq \frac{\tau_{1,S-1}^2}{m}\cdot\left(\frac{\tau_2 m}{\tau_{1,0}^2}\tilde{D}^0 + (1+2\beta(b))\,CLB_p(x^*, z_0)\right)\\ &= \frac{4}{m(S+3)^2}\Big(2m(F(\tilde{x}^0) - F(x^*)) + (1+2\beta(b))\,CLB_p(x^*, z_0)\Big).\end{aligned} \quad (A.35)$$

Meanwhile, by setting $\tilde{x}_0 = z_0 = 0$, using the optimality condition $0 \in \nabla f(x^*) + \partial\lambda\|x^*\|_1$ we have

$$\begin{aligned}F(\tilde{x}^0) - F(x^*) &= f(0) - f(x^*) - \lambda\|x^*\|_1\\ &\overset{\text{\textcircled{1}}}{\leq} \langle\nabla f(x^*), 0 - x^*\rangle + \frac{L}{2}\|x^*\|_1^2 - \lambda\|x^*\|_1\\ &\overset{\text{\textcircled{2}}}{=} \langle-\partial\lambda\|x^*\|_1, -x^*\rangle + \frac{L}{2}\|x^*\|_1^2 - \lambda\|x^*\|_1\\ &\overset{\text{\textcircled{3}}}{\leq} \|\partial\lambda\|x^*\|_1\|_\infty\|x^*\|_1 + \frac{L}{2}\|x^*\|_1^2 - \lambda\|x^*\|_1\\ &\overset{\text{\textcircled{4}}}{\leq} \lambda\|x^*\|_1 + \frac{L}{2}\|x^*\|_1^2 - \lambda\|x^*\|_1\\ &= \frac{L}{2}\|x^*\|_1^2,\end{aligned} \quad (A.36)$$

where \textcircled{1} is by the smoothness assumption of $f(x)$, \textcircled{2} is by selecting the subgradient of $\lambda\|x^*\|_1$ with $-\partial\lambda\|x^*\|_1 = \nabla f(x^*)$, \textcircled{3} is by lemma 12, \textcircled{4} is by using the property of subgradient $\|\partial\lambda\|x^*\|_1\|_\infty \leq \lambda$. In addition, for $1 < p \leq 2$,

$$\begin{aligned}B_p(x^*, z_0) &= B_p(x^*, 0) = \frac{1}{2}\|x^*\|_p^2 - \frac{1}{2}\|0\|_p^2 - \langle\nabla\frac{1}{2}\|0\|_p^2, x^* - 0\rangle\\ &= \frac{1}{2}\|x^*\|_p^2 \leq \frac{1}{2}\|x^*\|_1^2.\end{aligned} \quad (A.37)$$

Furthermore, minimizing $C = \frac{d^{\frac{2\delta}{1+\delta}}}{\delta}$ w.r.t $\delta$, we get $\delta = \log(d) - 1 - \sqrt{(\log(d)-1)^2 - 1}$ and $p = 1 + \delta = \log(d) - \sqrt{(\log(d)-1)^2 - 1} \in (1, 2]$,

Then combing (A.35), (A.36) and (A.37), we get the final result.

$$\mathbb{E}[F(\tilde{x}^S) - F(x^*)] \leq \frac{4}{(S+3)^2}\left(1 + \frac{1+2\beta(b)}{2m}C\right)L\|x\|_1^2. \quad (A.38)$$

$\square$